\documentclass[journal,twoside,web]{IEEEtran}

\usepackage{xcolor}
\usepackage{algorithmic}%
\usepackage{textcomp}

\usepackage[english]{babel}

\usepackage{amsmath, amssymb, amsfonts, amstext, array}

\usepackage{cases,bm,bbm}
\usepackage[mathscr]{eucal}
\usepackage{theorem,longtable,tabularx,booktabs}

\newtheorem{theorem}{Theorem}
\newtheorem{definition}{Definition}
\newtheorem{proposition}{Proposition}
\newtheorem{lemma}{Lemma}
\newtheorem{corollary}{Corollary}

{\theorembodyfont{\upshape}
\newtheorem{remark}{Remark}
\newtheorem{example}{Example}
}

\newcommand\mydots{\makebox[1em][c]{.\hfil.\hfil.}}

\def\qed{\hfill \vrule height 5pt width 5pt depth 0pt \medskip}
\newcommand{\proof}[1]{\noindent%
\ifx&#1&%
{\bf%
Proof. }
\else 
{\bf%
Proof of #1. }
\fi}

\newcommand{\beq}{\begin{equation}}
\newcommand{\eeq}{\end{equation}}
\newcommand{\beqa}{\begin{eqnarray}}
\newcommand{\eeqa}{\end{eqnarray}}
\newcommand{\beqan}{\begin{eqnarray*}}
	\newcommand{\eeqan}{\end{eqnarray*}}

\newcommand{\pde}[2]{ \frac{\partial #1}{\partial #2} }
\newcommand{\ppde}[2]{\frac{\partial^2 #1}{\partial #2^2}}
\newcommand{\ptde}[2]{\frac{\partial^3 #1}{\partial #2^3}}

\newcommand{\abs}[1]{ \left\lvert #1 \right\rvert }
\newcommand{\norm}[1]{ \| #1 \|}

\usepackage{dsfont,bbold}
\usepackage{graphicx} 
\usepackage{subfig}

\usepackage{enumitem}
\setlist[description]{leftmargin=\parindent}
\usepackage{xspace}

\newcommand{\qedh}{\tag*{\scriptsize$\blacksquare$}}

\DeclareMathOperator*{\argmax}{arg\,max}

\newcommand{\R}{\mathbb{R}}

\newcommand{\1}{\mathds{1}}
\newcommand{\0}{\mathbb{0}}

\newcommand{\diag}[1]{\text{diag}\left\{#1\right\}}
\newcommand{\sign}[1]{\text{sign}\left(#1\right)}
\newcommand{\range}[1]{\text{range}\,#1}
\newcommand{\card}[1]{ \text{card}(#1) }

\newcommand{\G}{\mathcal{G}}
\newcommand{\edgeSet}{\mathcal{E}}
\newcommand{\nodeSet}{\mathcal{V}}

\newcommand{\frustration}{{\epsilon(\G)}}

\newcommand{\spectrum}[1]{\Lambda(#1)}

\newcommand{\Lnorm}{\mathcal{L}}

\newcommand{\sym}[1]{{#1}_\text{sym}}
\newcommand{\dpsi}{{\Psi_{x}}}

\renewcommand{\ker}[1]{\text{ker}(#1)}
\renewcommand{\range}[1]{\text{range}(#1)}

\newcommand{\Sbest}{S_{\epsilon}}
\newcommand{\half}{\frac{1}{2}}
 
\newcommand{\xk}{{x_{k}}}
\newcommand{\xkp}{{x_{k+1}}}
\newcommand{\psik}{{\psi_k}}
\newcommand{\psikp}{{\psi_{k+1}}}
\newcommand{\psiki}{{\psi_{k,i}}}
\newcommand{\psikpi}{{\psi_{k+1,i}}}
\newcommand{\psiD}{{\psi_\Delta}}
\newcommand{\psiDi}{{\psi_{\Delta,i}}}
\newcommand{\der}{\text{d}}
\newcommand{\Lpi}{{L_{\pi}}}

\newcommand{\pid}{\pi_{1,d}}

\begin{document}
\newcolumntype{Z}{>{$}c<{$}}

\title{The role of frustration in collective decision-making dynamical processes on multiagent signed networks}

\author{Angela Fontan and Claudio Altafini
	\thanks{Work supported in part by a grant from the Swedish Research Council (grant n. 2015-04390). A preliminary version of this paper was presented at the 57th IEEE Conference on Decision and Control in 2018.}
	\thanks{A. Fontan and C. Altafini are with the Division of Automatic Control,
		Department of Electrical Engineering, Link\"{o}ping University, SE-58183 Link\"{o}ping,	Sweden, E-mail: $\{$angela.fontan, claudio.altafini$\}$@liu.se.} }

\maketitle

\begin{abstract}
	In this work we consider a collective decision-making process in a network of agents described by a nonlinear interconnected dynamical model with sigmoidal nonlinearities and signed interaction graph. 
	The decisions are encoded in the equilibria of the system.
	The aim is to investigate this multiagent system when the signed graph representing the community is not structurally balanced and in particular as we vary its frustration, i.e., its distance to structural balance. 
	The model exhibits bifurcations, and a ``social effort'' parameter, added to the model to represent the strength of the interactions between the agents, plays the role of bifurcation parameter in our analysis.	
	We show that, as the social effort increases, the decision-making dynamics exhibits a pitchfork bifurcation behavior where, from a deadlock situation of ``no decision'' (i.e., the origin is the only globally stable equilibrium point), two possible (alternative) decision states for the community are achieved (corresponding to two nonzero locally stable equilibria).
	The value of social effort for which the bifurcation is crossed (and a decision is reached) increases with the frustration of the signed network.
\end{abstract}

\section{Introduction}
\label{Introduction}
In this paper we want to study a nonlinear model for decision-making in a community of agents where antagonistic interactions may exist between the agents. 
Indeed, while collaboration between agents is often assumed in order to reach a common decision (for instance in applications such as collective behavior in animal groups \cite{GrayAl2018,SwainCouzinLeonard2012}, cooperative control in robotics \cite{RenBeard2008,RenCao2011}, or opinion forming \cite{HegselmannKrause2002,JiaAl2015}), there are applications in which restricting to collaborative interactions means oversimplifying the relationship among the agents \cite{WassermanFaust1994,EasleyKleinberg2010}. 
Classes of multiagent systems in which the presence of antagonism is plausible include for instance ``social networks'', i.e., groups of individuals interacting and exchanging opinions in a friendly/unfriendly manner or trusting/mistrusting each other.
Other scenarios in which antagonism is unavoidable are team games, where different teams have to compete against each other, or parliamentary democracies, where parties can be allied or rival.

Signed networks \cite{Harary1953,Zaslavsky1982} are a natural framework to model a community of agents where both cooperative and antagonistic interactions coexist: a positive sign labeling an edge between two agents represents a friendly (or cooperative) relationship, while a negative sign labeling an edge an unfriendly (or competitive) relationship.
If the group of agents can be divided into two subgroups such that the agents inside each group are mutual friends (i.e., they are linked by edges with positive weight) while the agents across the two subgroups are enemies (i.e., they are linked by edges with negative weight), we say that the network is structurally balanced \cite{CartwrightHarary1956,Altafini2013TAC}.
If we assume that a network is undirected and connected, an equivalent condition to structural balance is that the smallest eigenvalue of the normalized signed Laplacian $\Lnorm$ is zero, $\lambda_1(\Lnorm)=0$ in the notation we introduce below.
As for instance the works \cite{FacchettiIaconoAltafini2011,KunegisAl2010,Kunegis2014,FontanAltafini2021} show, real signed social networks are in general not structurally balanced.

To model the evolution of the opinions of the agents in a community represented as a signed social network we use the model of opinion forming previously introduced in \cite{GrayAl2018,AbaraTicozziAltafini2018,FontanAltafini2018}. 
This model is characterized by sigmoidal and saturated nonlinearities, describing how the agents transmit their opinion to their neighbors. It has a (signed) Laplacian-like structure at the origin and it is endowed with a social effort parameter $\pi$ which in our analysis plays the role of bifurcation parameter. 
Our aim is to study how the strength of the commitment among the agents, represented by $\pi$, affects the presence and stability of the equilibrium points of the system, which represent the decision states for the community.
Under our assumptions, the system is monotone \cite{Smith1988} if and only if the corresponding signed social network is structurally balanced. 
In this case the behavior of the system can be easily deduced from \cite{GrayAl2018,AbaraTicozziAltafini2018,FontanAltafini2018}, where the authors consider a cooperative system (i.e., only friendly interactions exist between the agents), which is a particular case of monotone system. 
In this case the analysis shows that for increasing values of the social effort parameter $\pi$, the system undergoes two sequential pitchfork bifurcations: after the first bifurcation the number of equilibria jumps from one to three, while after the second bifurcation multiple (more than three) equilibrium points arise.
In particular when crossing the first bifurcation the system passes from having the origin as globally asymptotically stable equilibrium to a situation in which two nonzero locally stable equilibria exist while the origin becomes a saddle point.
This situation is maintained up to the second bifurcation where novel equilibria, stable or unstable, appear.
In the context of social interactions this behavior can be interpreted as follows: if the social effort between the agents is small then no decision is achieved (the origin is the only attractor), while two alternative decision states can be reached if the agents have the ``right'' amount of commitment.
However, by further increasing the social effort, the agents may fall in a situation of overcommitment where multiple (more than 2) decisions are possible.
For cooperative networks the first threshold value is fixed and constant, while the second threshold value depends on the algebraic connectivity of the network. 

We show in this work that if we consider signed networks that are structurally unbalanced then, while the qualitative behavior of the system does not change, the value of social effort parameter for which the first bifurcation is crossed is no longer constant but grows with the smallest eigenvalue of the normalized signed Laplacian of the network, which for structurally unbalanced networks is strictly positive ($\lambda_1(\Lnorm)>0$).
In particular, 
its value increases with the amount of ``frustration'' encoded in the signed network, i.e., with the amount of ``disorder'' that the negative edges introduce in a network, see \cite{FacchettiIaconoAltafini2011} for a more thorough statistical physics interpretation. 
First introduced by Harary \cite{Harary1953,Harary1959} and denoted ``line index of balance'', the frustration is a standard measure to express the distance of a signed network from a structurally balanced state and is defined as the minimum (weighted) sum of the negative edges that need to be removed in order to obtain a structurally balanced network, see \cite{FontanAltafini2018CDC} for details.

For our model of decision-making, this means that when we consider signed networks with higher frustration, the first bifurcation is crossed at higher values of the social effort parameter $\pi$, meaning that a higher commitment is required from the agents in order to converge to a nontrivial equilibrium point.
From a sociological point of view, the result admits a fairly reasonable interpretation: the more in a community there are ``unresolved tensions'' among the agents (i.e., unbalanced interactions, as measured by the frustration), the more commitment is required by the agents to achieve a common nontrivial decision and to ``escape'' the (trivial) zero equilibrium point.
On the other hand, the value of social effort for which the second bifurcation is crossed is independent of the frustration of the network \cite{FontanAltafini2018CDC}, meaning that for highly frustrated graphs the range of social commitment values for which only two nontrivial equilibria are present shrinks. 
As a concrete application of these results, in a recent work (see \cite{FontanAltafini2021}) we have described the process of government formation in parliamentary democracies as a collective decision-making process where the members of the parliament (the agents) are required to cast a vote of confidence (the decision) to a candidate cabinet coalition.
In this context, the social effort parameter $\pi$ is a proxy for the complexity of the government negotiations (measured as duration of the negotiation phase), while a signed network describes the composition of the parliament after each election with signs representing party alliances/rivalries.
These ``parliamentary networks'' are in general not structurally balanced and their frustration correlates well with the duration of the government negotiation processes.

Because of the nonlinearities the behavior of our system is fundamentally different from that of \cite{Altafini2013TAC}. In the case of \cite{Altafini2013TAC} in fact, structural balance leads to bipartite consensus and structural unbalance to asymptotic stability. In our case, instead, balanced and unbalanced cases are qualitatively similar, with only the bifurcation point gradually moving to higher values of social commitment as the frustration grows.
In this respect, the model we present here has a more reasonable behavior than the one in \cite{Altafini2013TAC}, at the cost of a higher complexity.

Even though the behavior of the system in the structurally unbalanced case is qualitatively similar to the structurally balanced case, the technical tools that must be used to show the results become much more challenging because the system is no longer monotone.
An important technical contribution of this paper is in fact to develop methods able to perform a global state space analysis of a broad class of nonlinear nonmonotone interconnected systems which are not diagonally dominant.
Familiar examples of Hopfield-like neural networks fall in this category \cite{Hopfield1984,ZhangWangLiu2014}.
Another noteworthy result we obtain is a description of the region in which all equilibria of the system must be contained. In particular also the upper bound to the 1-norm of the equilibria we provide depends on the frustration of the signed network.

The paper investigates also a discrete-time version of our multiagent decision-making system. Such extension is nontrivial in several directions: for instance new phenomena, like period-2 limit cycles, appear in the discrete-time case. Also the techniques that must be used to prove the results are largely different from those of the continuous-time case. 
In particular, we show that the first bifurcation occurring at the origin is either a pitchfork or a period-doubling bifurcation, depending on the relative positions of the corresponding threshold values for the social effort parameter $\pi$, where the value for which a pitchfork bifurcation is crossed is the same as in the continuous-time case.
Interestingly, we show that if the signed network has zero or small frustration, the value of $\pi$ for which a period-doubling bifurcation is crossed is always bigger than the usual bifurcation threshold.

The rest of the paper is organized as follows: in Section~\ref{sec:Preliminaries} we introduce preliminary material. In Sections~\ref{DecisionMakingSystem} and \ref{sec:DT} we present our results for collective decision-making over signed networks in continuous- and discrete-time, respectively. 
The results are discussed and interpreted in Section~\ref{sec:Results}.
Numerical simulations and examples are shown in Section~\ref{Examples}. 
Technical preliminaries (useful Lemmas and Theorems) and most of the proofs are put in the Appendices at the end of the paper.

A preliminary version of this work appears in the conference proceedings of CDC 2018 \cite{FontanAltafini2018CDC}.
The new contributions are a necessary condition for the existence of nontrivial equilibria, the proof that these equilibria are locally asymptotically stable, and the description of the region in which all equilibria must be contained. All the material on the discrete-time version of our decision-making model is presented here for the first time.

\section{Preliminaries}
\label{sec:Preliminaries}

\subsection{Notation and linear algebra}
Given a matrix $A= [a_{ij}] \in \R^{n\times n}$, $A\ge 0$ means element-wise nonnegative, i.e., $a_{ij}\ge 0$ for all $i,j=1,\mydots,n$, while $A>0$ means element-wise positive, i.e., $a_{ij}>0$ for all $i,j=1,\mydots,n$.
The spectrum of $A$ is denoted $\spectrum{A}= \{\lambda_1(A),\dots,\lambda_n(A)\}$, where $\lambda_i(A)$, $i=1,\mydots,n$, are the eigenvalues of $A$.
A matrix $A$ is called irreducible if there does not exist a permutation matrix $P$ s.t. $P^T A 
P$ is block triangular.
If $x,y\in \R^n$ then $x\ge y$ ($x>y$) means that $x_i \ge y_i$ (resp., $x_i>y_i$) for all $i=1,\mydots,n$.
Given two matrices $A,B\in \R^{n\times n}$, the notation $A\sim B$ means that $A$ and $B$ are similar, and hence that they have the same eigenvalues.
Given a diagonal positive definite matrix $D$, we denote the unique (diagonal) positive definite square root of $D$ by $D^{\half}$.
The symbol $\1$ indicates the vector of ones ($\1_m$ is used when the dimension $m$ is not clear from the context) and $0_{n,m}$ the $n\times m$ zero matrix ($0$ if it is clear from the context).

\subsection{Signed graphs}
\label{SignedGraphs}
Let $\G=(\nodeSet,\edgeSet)$ be a graph with vertex set $\nodeSet$ (such that $\card{\nodeSet}=n$) and edge set $\edgeSet\subseteq \nodeSet\times \nodeSet$. Let $A=[a_{ij}]\in \R^{n\times n}$ be the adjacency matrix of $\G$, i.e., $a_{ij}\ne 0$ if and only if $(j,i)\in \edgeSet$.
In this work we consider undirected and connected graphs without self-loops. 

A graph $\G$ is \textit{signed} if each of its edges is labeled by a sign, that is, $\sign{a_{ij}}=\sign{a_{ji}}=\pm 1$ if $(i,j)\in \edgeSet$.
The \textit{signed Laplacian} of a graph $\G$ is the symmetric matrix $L= \Delta-A$, where $\Delta = \diag{\delta_1,\dots,\delta_n}$ and each diagonal element $\delta_i$ is given by $\delta_i= \sum_{j=1}^{n} \abs{a_{ij}}$, $i=1,\mydots,n$ \cite{Altafini2013TAC}.
The \textit{normalized signed Laplacian} of a graph $\G$, see \cite{LiLi2009,HouLiPan2003}, is the non-symmetric (symmetrizable, see Appendix~\ref{sec:TechnicalPreliminaries} for a definition) matrix defined as 
\begin{equation}
\mathcal L = \Delta^{-1}L = I- \Delta^{-1}A.
\label{eqn:NormalizedLaplacian}
\end{equation}
Notice that since the graph $\G$ is connected, it does not have isolated vertices, hence $\delta_i\ne0$ for all $i$ and the matrix $\Delta^{-1}$ is well-defined and positive definite.

All the matrices we consider in this work are either symmetric (e.g., $A$ and $L$) or symmetrizable (e.g., $\Lnorm$), hence they have real eigenvalues, which we assume to be arranged in a nondecreasing order. Let $\lambda_i(A)$, $\lambda_i(L)$ and $\lambda_i(\Lnorm)$, $i=1,\mydots,n$, be the eigenvalues of $A$, $L$ and $\Lnorm$, respectively.
By construction the eigenvalues of the signed Laplacian $L$ and the normalized signed Laplacian $\Lnorm$ are nonnegative, which can be easily shown using the Ger\v{s}gorin's Theorem, see \cite[Thm 6.1.1]{HornJohnson2013}.

A cycle of a signed graph $\G$ is said \textit{positive} if it contains an even number of negative edges, \textit{negative} otherwise.
A graph $\G$ is \textit{structurally balanced} if all its cycles are positive.
Equivalent conditions for $\G$ (connected) to be structurally balanced are the following \cite{Altafini2013TAC}: (i) there exists a partition of the node set $\nodeSet=\nodeSet_1\cup \nodeSet_2$ such that every edge between $\nodeSet_1$ and $\nodeSet_2$ is negative and every edge within $\nodeSet_1$ or $\nodeSet_2$ is positive; (ii) there exists a signature matrix $S=\diag{s_1,\mydots,s_n}$ with diagonal entries $s_i=\pm 1$ ($i=1,\mydots,n$), such that $S\Lnorm S$ has all nonpositive off-diagonal entries; (iii) $\lambda_1(\Lnorm)=0$.
The \textit{frustration index} of a signed graph $\G$ is defined as
\begin{equation}
\frustration
= \min_{\substack{S=\diag{s_1,\mydots,s_n}\\s_i=\pm 1\; \forall i}}\;\;
\frac{1}{2} \sum_{i\ne j} \,[\,\abs{\Lnorm}+S\Lnorm S\,]_{ij},
\label{eqn:frustration}
\end{equation}
where $[\cdot]_{ij}$ indicates the $i,j$ element and $\abs{\cdot}$ the element-wise absolute value, and it provides a measure of the distance of $\G$ from a structurally balanced state \cite{FontanAltafini2018CDC}. If $\G$ is structurally balanced, $\frustration = 0$.

\subsection{Monotone systems}
Consider the system 
\begin{equation}
\dot x = f(x),\quad x(0)=x_0
\label{eqn:SystemExample}
\end{equation}
where $f$ is a continuously differentiable function defined on a convex open set $\mathcal U \subseteq \R^n$. Let $\varphi(t,\bar{x})$ be the solution $x(t)$ of \eqref{eqn:SystemExample} s.t. $x(0)= \bar{x}$.

Let $S$ be a signature matrix, i.e., $S= \diag{s_1,\mydots,s_n}$ with $s_i=\pm1$ $\forall i$, and let $S\R^n$ indicate an orthant of $\R^n$, $S\R^n = \{x\in \R^n: s_i x_i \ge 0,\,i=1,\mydots,n\}$. The partial ordering $\le_S$ is preserved by the solution operator $\varphi(t,\cdot)$ and the system~\eqref{eqn:SystemExample} is type $S\R^n$ monotone if whenever $\bar{x},\bar{y}\in \mathcal{U}$ with $\bar{x}\le_S \bar{y}$ then $\varphi(t,\bar{x}) \le_S \varphi(t,\bar{y})$ for all $t\ge 0$ \cite{Smith1988}.
\begin{lemma}[2.1 in \cite{Smith1988}]
	\label{lemma:monotone}
	If $f\in C^1(\mathcal{U})$ where $\mathcal U$ is open and convex in $\R^n$ then $\varphi(t,\cdot)$ preserves the partial ordering $\le_S$ for $t\ge 0$ if and only if $S \pde{f}{x}(x) S$ has nonnegative off-diagonal elements for every $x\in \mathcal U$.
\end{lemma}
Therefore, a system~\eqref{eqn:SystemExample} is monotone if and only if the graph described by the Jacobian $\pde{f}{x}$ as adjacency matrix is structurally balanced with fixed $S$ $\forall x\in \mathcal{U}$.

\section{Decision-making in antagonistic multiagent systems in continuous-time}
\label{DecisionMakingSystem}

\subsection{Problem formulation}
\label{ProblemFormulation}
To model the process of decision-making in a community of $n$ agents represented by a signed network $\G$, we consider the following class of nonlinear interconnected systems,
\begin{equation}
\dot x = -\Delta x+\pi A \psi(x),\quad x\in \R^n.
\label{eqn:System_A}
\end{equation}
The state vector $x= [x_1\,\cdots\,x_n]^T\in \R^n$ represents the agents' opinions, $A=[a_{ij}]$ is the adjacency matrix of the network $\G$ and describes how the agents interact with each other, $\Delta=\diag{\delta_1,\dots,\delta_n}$, $\pi>0$ is a positive scalar parameter and $\psi(x)=[\psi_1(x_1)\,\cdots\, \psi_n(x_n)]^T$.
Each nonlinear function $\psi_i(x_i)$ describes how an agent $i$ transmits its opinion $x_i$ to its neighbors in the network.
This term is then weighted first by the element $a_{ij}$, describing the influence between agents $i$ and $j$ (positive/friendly if $a_{ij}>0$ or negative/unfriendly if $a_{ij}<0$), and then by the parameter $\pi$ representing the global ``social effort'' or ``strength of commitment'' among the agents \cite{GrayAl2018}.
The equilibria of the system represent the decision states for the community.

We assume that the signed network $\G$ is undirected (two agents able to influence each other's opinion share the same amount of trust/distrust in each other), connected (there are no isolated agents) and without self-loops, meaning that the signed adjacency matrix $A$ is symmetric, irreducible and with null diagonal. 
We also assume that a Laplacian-like assumption relates $\Delta$ and $A$, $\delta_i=\sum_{j} \abs{a_{ij}}$.
Finally, we assume that each nonlinear function $\psi_i(x_i):\R\to \R$ of the vector $\psi(x)$ satisfies the following conditions
\begin{gather}
\psi_i(x_i)=-\psi_i(-x_i),\,\forall x_i\in \R\quad\text{(odd)}
\tag{A.1}
\label{assumption:1psiOdd}
\\
\pde{\psi_i}{x_i}(x_i)>0,\,\forall x_i\in \R\;\text{and }\pde{\psi_i}{x_i}(0)=1\quad\text{(monotone)}
\tag{A.2}
\label{assumption:2psiMonotone}
\\
\lim_{x_i\to\pm \infty} \psi_i(x_i)=\pm 1\quad\text{(saturated)}
\tag{A.3}
\label{assumption:3psiSaturated}
\\
\psi_i(x_i) \; 
\begin{cases}
\text{strictly convex} & \forall\, x_i<0\\
\text{strictly concave}& \forall\, x_i>0
\end{cases} \quad \text{(sigmoidal)}.
\tag{A.4}
\label{assumption:4Sigmoidal}
\end{gather}
The system~\eqref{eqn:System_A} can be rewritten in a ``normalized'' form,
\begin{equation}
\dot x = \Delta\left[-x+\pi H \psi(x)\right],\quad x \in \R^n,
\label{eqn:System_CT}
\end{equation}
where we consider the normalized interaction matrix $ H := \Delta^{-1} A $.
The Jacobian of \eqref{eqn:System_CT} is $J(x) =-\Delta(I - \pi H \pde{\psi}{x}(x))$ which at the origin for $\pi=1$ reduces to $J=-\Delta \,\Lnorm$ (where $\Lnorm$ is the normalized signed Laplacian of the network); hence, under our assumptions and from Lemma~\ref{lemma:monotone}, the system~\eqref{eqn:System_CT} is monotone if and only if the signed network $\G$ is structurally balanced.

Our aim is use bifurcation analysis to investigate how the social effort parameter $\pi$ (our bifurcation parameter) affects the presence of the equilibrium points of the system~\eqref{eqn:System_CT}.

\subsection{Structurally balanced case}
Previous works, such as \cite{AbaraTicozziAltafini2018,FontanAltafini2018}, have studied the behaviour of the system~\eqref{eqn:System_CT} when the adjacency matrix $A$ of the network is nonnegative, i.e., when the system is cooperative \cite{Smith1988}. These results, summarized in the following theorem, still hold when the system is in general monotone, that is, when the network $\G$ described by the matrix $A$ is structurally balanced.
\begin{theorem}[\cite{FontanAltafini2018}]
	\label{thm:CT_SB_summary}	
	Consider the system~\eqref{eqn:System_CT} where each nonlinear function $\psi_i(\cdot)$, $i=1,\mydots,n$, satisfies the properties \eqref{assumption:1psiOdd}$\div$\eqref{assumption:4Sigmoidal}. 
	Assume that the signed graph $\G$ is structurally balanced and let $S$ be the signature matrix s.t. $S\Lnorm S$ has all nonpositive off-diagonal entries ($\abs{A} = S A S$).
	\begin{enumerate}[label=(\roman*)]
		\item When $\pi<1$ the origin is the unique equilibrium point and it is asymptotically stable.
		\item When $\pi=\pi_1=1$ the system undergoes a pitchfork bifurcation, the origin becomes unstable and two new equilibria appear, in the orthants described by $S$ and $-S$, respectively, denoted $S\R^n_+$ and $S\R^n_-$. These equilibria are locally asymptotically stable with domain of attraction at least equal to $S\R^n_+$ and $S\R^n_-$, respectively.
		\item If $\lambda_2(\Lnorm)<1$ and simple, when $\pi=\pi_2=\frac{1}{1-\lambda_2(\Lnorm)}$ the system undergoes a second pitchfork bifurcation, and new equilibria in other orthants of $\R^n$ appear, which may be stable or unstable.
	\end{enumerate}
\end{theorem}

\subsection{Structurally unbalanced case}
In this section we want to introduce our novel results, i.e., the extension of Theorem~\ref{thm:CT_SB_summary} to signed networks which are structurally unbalanced: we show that, by redefining the threshold values $\pi_1$ and $\pi_2$, the system~\eqref{eqn:System_CT} behaves in a similar manner as the one described in Theorem~\ref{thm:CT_SB_summary}.

Theorem~\ref{thm:CT_SUB_summary} summarizes our findings.
We proceed as follows: first, in (i), we prove that the origin is the unique equilibrium point for the system when $\pi<\pi_1$ and it is globally asymptotically stable, where $\pi_1$ depends on the smallest eigenvalue of the normalized signed Laplacian $\Lnorm$. 
Then, in (ii) we show that when $\pi = \pi_1$ the system undergoes a pitchfork bifurcation and two new equilibria appears, which are locally asymptotically stable for all values of the bifurcation parameter in the interval $(\pi_1,\pi_2)$, where $\pi_2$ depends on the second smallest eigenvalue of the normalized signed Laplacian $\Lnorm$.
Similarly to \cite{GrayAl2018,FontanAltafini2018}, the proof relies on bifurcation theory.
Lack of monotonicity however implies that most of the proofs require different arguments than those used in \cite{GrayAl2018,FontanAltafini2018}.
At $\pi_2$ the system bifurcates again and new equilibria appear, see (iii).
\begin{theorem}\label{thm:CT_SUB_summary}
Consider the system~\eqref{eqn:System_CT} where each nonlinear function $\psi_i(\cdot)$, $i=1,\mydots,n$, satisfies the properties \eqref{assumption:1psiOdd}$\div$\eqref{assumption:4Sigmoidal}.
Assume that the signed graph $\G$ is structurally unbalanced with normalized signed Laplacian $\Lnorm$. 	
\begin{enumerate}[label=(\roman*)]
	\item When $\pi\le \pi_1=\frac{1}{1-\lambda_1(\Lnorm)}$, the origin is the unique equilibrium point of \eqref{eqn:System_CT} and it is globally asymptotically stable.

	\item Let $\lambda_1(\Lnorm)$ be simple,
	\begin{enumerate}
		\item[(ii.1)] (existence): when $\pi$ crosses $\pi_1$ the system undergoes a pitchfork bifurcation and two new equilibria ($x^*$ and $-x^*$) appear;
			
		\item[(ii.2)] (stability): when $\pi>\pi_1$ the origin is an unstable equilibrium point, while the equilibria $\pm x^\ast\ne 0$ are locally asymptotically stable for all values of $\pi\in (\pi_1,\pi_2)$, with $\pi_2 =\frac{1}{1-\lambda_{2}(\Lnorm)}$;
		
		\item[(ii.3)] (uniqueness): when $\pi\in (\pi_1,\pi_2)$, the system admits exactly three equilibria, the origin and the two nontrivial equilibrium points $\pm x^\ast\ne 0$.
	\end{enumerate}
	
	\item If $\lambda_2(\Lnorm)$ is simple, when $\pi=\pi_2$ the system undergoes a second pitchfork bifurcation and new equilibria appear.
\end{enumerate}
\end{theorem}
Proof in Appendix~\ref{proof:CT_SUB_Summary}.
\begin{remark}\label{remark:Eq}
It follows from the assumption \eqref{assumption:1psiOdd} that if the system~\eqref{eqn:System_CT} admits an equilibrium point $x^*\ne0$, then $-x^*$ is also an equilibrium point.
\end{remark}
\begin{remark}
Differently from Theorem~\ref{thm:CT_SB_summary}(iii), in Theorem~\ref{thm:CT_SUB_summary}(iii) the assumption $\lambda_{2}(\Lnorm)<1$ is not needed: if the network $\G$ is structurally unbalanced and connected it is always true that $\lambda_2(\Lnorm)<1$, as shown in Lemma~\ref{lemma:lambda2} below.
Therefore, if $\lambda_2(\Lnorm)$ is simple, $\pi_2 = \frac{1}{1-\lambda_2(\Lnorm)}$ is always well-defined (i.e., strictly positive and greater than $\pi_1$).
On the other hand, examples of structurally balanced graphs for which $\lambda_2(\Lnorm)>1$ are complete graphs, whose adjacency matrix is a Euclidean distance matrix. 
This means that in the structurally balanced case when $\lambda_2(\Lnorm)>1$ the system~\eqref{eqn:System_CT} admits only $3$ equilibrium points ($0$, $\pm x^\ast$) for all values of $\pi>\pi_1$ and that the trajectories converge either to $x^\ast$ or $-x^\ast$.
However, this situation can never happen in the structurally unbalanced case: if $\pi$ is ``large enough'' (i.e., it is above the threshold $\pi_2$) the system~\eqref{eqn:System_CT} will always admit new equilibria (other than $0$, $\pm x^\ast$), which may be attractors.
\end{remark}
\begin{lemma}\label{lemma:lambda2}
Let $\G$ be a signed connected network with normalized signed Laplacian $\Lnorm$.
If $\G$ is structurally unbalanced, $\lambda_2(\Lnorm)<1$.
\end{lemma}
Proof in Appendix~\ref{proof:lambda2}.

To conclude this part, we show that for $\pi>\pi_1$ the 1-norm of the equilibria of the system is upper bounded by $\pi (n-2\frustration)$, where $\frustration$ is the frustration of the signed network.
Moreover, if the matrix $\Lnorm$ is symmetric (i.e., if $\Delta = \delta I$), we can show that the solutions of \eqref{eqn:System_CT} are all bounded and converge to a set $\Omega_{\frustration}$, which implies that all the equilibria of the system~\eqref{eqn:System_CT} belong to $\Omega_{\frustration}$. 
\begin{theorem}\label{thm:CT_BoundedEq}	
Consider the system~\eqref{eqn:System_CT} where each nonlinear function $\psi_i(\cdot)$, $i=1,\mydots,n$, satisfies the properties \eqref{assumption:1psiOdd}$\div$\eqref{assumption:4Sigmoidal}. 
Let $\frustration$ be the frustration of the signed network $\G$ defined in \eqref{eqn:frustration}, and $\Lnorm$ its normalized signed Laplacian.
\begin{enumerate}[label=(\roman*)]
	\item If $x^*$ is an equilibrium point of \eqref{eqn:System_CT}, then $\norm{x^*}_1\le\pi(n-2\frustration)$.
		
	\item Let $\pi>\pi_1$. Under the assumption that $\Lnorm$ is symmetric (i.e., $\Delta = \delta I$), the trajectories of \eqref{eqn:System_CT} asymptotically converge to the set $\Omega_{\frustration}$, where
	\begin{equation*}
		\Omega_{\frustration}= \{x\in \R^n: \norm{x}_1\le \pi (n-2\frustration)\}.
	\end{equation*}
\end{enumerate}
\end{theorem}
Proof in Appendix~\ref{sec:BoundedTraj}.
\begin{proposition}\label{prop:fpi1pi2}
Let $\G$ be a signed graph with normalized signed Laplacian $\Lnorm$, and assume that $\Lnorm$ is symmetric (i.e., $\Delta = \delta I$).
Then it is possible to derive an upper bound for the social effort at the first bifurcation point, $\pi_1$, which depends on the frustration of the network $\frustration$:
\begin{equation}
	1 \le \pi_1 \le \min\Bigl\{ \frac{n}{n-2\frustration},\pi_2\Bigr\}.
	\label{eqn:fpi1pi2}
\end{equation}
\end{proposition}
\noindent%
Proof in Appendix~\ref{sec:prop:fpi1pi2}. Notice that if the frustration is zero (i.e., the network is structurally balanced) then $\pi_1 = 1$.

\section{Discrete-time}\label{sec:DT}%
The Euler approximation of system~\eqref{eqn:System_CT} with step $\varepsilon$ is 
\begin{multline}
x_i(k+1) 
= (1-\varepsilon \delta_i)\, x_i(k)+\varepsilon\pi \sum_{j\ne i} a_{ij} \psi_j(x_j(k)),\\ i = 1,\mydots,n.
\label{eqn:System_DT_i}
\end{multline}
Let $\xk:=x(k) =[x_1(k) \, \cdots \, x_n(k)]^T$ and $\psi(\xk) = [\psi_1(x_1(k))\, \cdots \, \psi_n(x_n(k))]^T$.
Equation \eqref{eqn:System_DT_i} can be rewritten in a more compact form as follows:
\begin{equation}
\xkp = (I-\varepsilon \Delta) \xk+ \varepsilon \pi A \psi(\xk).
\label{eqn:System_DT}
\end{equation}
The Jacobian at the origin is given by:
\begin{equation}
J_\pi 
= I-\varepsilon \Delta+ \varepsilon \pi  A
= I-\varepsilon \Lpi
\label{eqn:Jpi}
\end{equation}
where
\begin{equation}
\Lpi = \Delta -\pi A = \Delta\bigl(I-\pi (I-\Lnorm)\bigr)
\label{eqn:Lpi}
\end{equation}
and $\Lnorm$ is the normalized signed Laplacian of the network.
As in Section~\ref{DecisionMakingSystem}, we want to study how the social effort parameter $\pi$ affects the existence of the equilibria of the system~\eqref{eqn:System_DT}, again relying on tools from bifurcation analysis \cite{Kuznetsov1998}.

A local bifurcation occurs at the origin if the Jacobian $J_\pi$ has an eigenvalue with absolute value equal to $1$ (that is, equal to $\pm 1$ since $J_\pi$ is symmetric and has real eigenvalues).
When $\pi$ is small and in particular is such that all the eigenvalues of $J_\pi$ have magnitude less than one, following the proof of \cite[Thm~2]{Zhaoetal2002} it is possible to prove (under the additional condition $\varepsilon \max_i \delta_i < 1$) that the origin is globally asymptotically stable and hence the unique equilibrium point of the system~\eqref{eqn:System_DT}.
As $\pi$ grows, the magnitude of the eigenvalues of $J_\pi$ increases and for values of $\pi$ such that $J_\pi$ has a simple eigenvalue $\lambda$ at $\pm 1$ the system~\eqref{eqn:System_DT} can undergo either a pitchfork ($\lambda=+1$) or a period-doubling ($\lambda=-1$) bifurcation \cite{Kuznetsov1998}.

Let $\lambda_i(J_\pi)$ and $\lambda_i(\Lpi)$, $i=1,\mydots,n$, be the eigenvalues of $J_\pi$ and $\Lpi$, respectively, which we assume to be arranged in a nondecreasing order. We denote $\pi_1$ the value of social effort for which the biggest eigenvalue of $J_\pi$ crosses $+1$ and $\pid$ the value of social effort for which the smallest eigenvalue of $J_\pi$ crosses $-1$:
\begin{equation}
\pi_1: \,\lambda_n(J_{\pi_1})=1,\quad \pid:\,\lambda_1(J_{\pid})=-1.
\label{eqn:pi1_pi1d}
\end{equation}
\begin{remark}\label{remark:eig_Lpi}
From \eqref{eqn:Jpi}, the biggest and smallest eigenvalues of $J_\pi $ are given by
\begin{align*}
\lambda_n(J_\pi ) = 1-\varepsilon \lambda_{1}(\Lpi),\quad 
\lambda_1(J_\pi) = 1-\varepsilon \lambda_{n}(\Lpi),
\end{align*}
which means that $J_\pi$ is Schur stable (i.e., its eigenvalues have magnitude strictly less than one) if and only if the following two conditions hold:
\begin{enumerate}[label=(\roman*)]
	\item $\lambda_{1}(\Lpi)>0$, that is, $\Lpi$ is positive definite. \\From \eqref{eqn:Lpi} and since $\Delta$ is positive definite, this condition is equivalent to $I-\pi (I-\Lnorm)$ having (strictly) positive eigenvalues, i.e., $0<1-\pi (1-\lambda_1(\Lnorm))$;
	
	\item $\varepsilon \lambda_{n}(\Lpi)<2$, that is, $\varepsilon \Lpi - 2 I$ is negative definite.
\end{enumerate}
Hence $\pi_1$ is the value of social effort for which the smallest eigenvalue of $\Lpi$ crosses $0$, while $\pid$ is the value of social effort for which the biggest eigenvalue of $\varepsilon\Lpi$ crosses $2$:
\begin{align}
\pi_1&: \,\lambda_1(L_{\pi_1})=0 \;\Rightarrow\;\pi_1 = \frac{1}{1-\lambda_1(\Lnorm)}
\label{eqn:pi1}
\\
\pid&:\,\lambda_n(L_{\pid})= \frac{2}{\varepsilon}.
\label{eqn:pi1d}
\end{align}
\end{remark}
In the analysis of the discrete-time model~\eqref{eqn:System_DT} (see Theorem~\ref{thm:DT_summary} below) it is relevant to know where $\pid$ lies compared with $\pi_1$. 
The next proposition shows that $\pi_1<\pid$ always holds if the network is structurally balanced ($\lambda_1(\Lnorm)=0$) or if it is structurally unbalanced but $\lambda_1(\Lnorm)>0$ is small. 
\begin{proposition}\label{prop:fpi1pi1d}
Assume that $\varepsilon \max_i \delta_i <1$. If any of the two following conditions on the signed graph $\G$ with normalized signed Laplacian $\Lnorm$ is satisfied:
\begin{enumerate}[label=(\roman*)]
	\item $\G$ is structurally balanced (i.e., $\lambda_1(\Lnorm)=0$), or
	\item $\G$ is structurally unbalanced and $\lambda_1(\Lnorm)<2-\lambda_n(\Lnorm)$,
\end{enumerate}
then $\pi_1<\pid$.
\end{proposition}
Proof in Appendix~\ref{sec:DT_fpi1pi1d}.

The next two lemmas show that if the system~\eqref{eqn:System_DT} admits a nontrivial equilibrium point then $\pi>\pi_1$ (Lemma~\ref{lemma:DT_Origin}), while if it admits a period-2 orbit then $\pi>\pid$ (Lemma~\ref{lemma:DT_Origin_Oscill}). 
\begin{lemma}
	\label{lemma:DT_Origin}
	Consider the system~\eqref{eqn:System_DT} where each nonlinear function $\psi_i(\cdot)$, $i=1,\mydots,n$, satisfies the properties \eqref{assumption:1psiOdd}$\div$\eqref{assumption:4Sigmoidal}. 
	If $x^\ast\ne 0$ is an equilibrium point of the system~\eqref{eqn:System_DT} then $\pi>\pi_1$.
\end{lemma}
Proof in Appendix~\ref{sec:DT_NecCondNontrivialEq}.
\begin{lemma}
	\label{lemma:DT_Origin_Oscill}
	Consider the system~\eqref{eqn:System_DT} where each nonlinear function $\psi_i(\cdot)$, $i=1,\mydots,n$, satisfies the properties \eqref{assumption:1psiOdd}$\div$\eqref{assumption:4Sigmoidal}. 
	If $\varepsilon <\frac{2}{\max_i \delta_i}$ and the system~\eqref{eqn:System_DT} admits a period-2 limit cycle ($\exists\, K>0$ such that $x_{k+2}=\xk$ for all $k\ge K$) then $\pi>\pid$.
\end{lemma}
Proof in Appendix~\ref{sec:DT_NecCondNontrivialEq}.
\begin{remark}
	The condition $\varepsilon <\frac{2}{\max_i \delta_i}$ imposed by Lemma~\ref{lemma:DT_Origin_Oscill} represents an upper bound on the step size $\varepsilon$ in the Euler approximation.
\end{remark}

We are now ready to state our results for the discrete-time system~\eqref{eqn:System_DT}, summarized in Theorem~\ref{thm:DT_summary}.
Similarly to Theorem~\ref{thm:CT_SUB_summary}, we show first that the origin is the unique equilibrium point for the system when $\pi<\min\{\pi_1,\pid\}$ and that it is globally asymptotically stable, see Theorem~\ref{thm:DT_summary}(i). 
However, differently from the continuous-time case, when $\pi$ crosses $\min\{\pi_1,\pid\}$ two different behaviors can happen.
If $\pi_1<\pid$ we expect the system~\eqref{eqn:System_DT} to undergo a pitchfork bifurcation when $\pi=\pi_1$ while if $\pi_1>\pid$ we expect a period-doubling bifurcation when $\pi=\pid$, see Theorem~\ref{thm:DT_summary}(ii).
The special case where $\pi_1=\pid$ is here not discussed, but the intuition is that a Neimark-Sacker bifurcation occurs at the origin when $\pi = \pid = \pi_1$ \cite{Kuznetsov1998}.
Observe also that $\pi>\pi_1$ and $\pi>\pid$ are necessary conditions (not only sufficient) in order for the system~\eqref{eqn:System_DT} to admit a nontrivial equilibrium or a periodic solution, respectively, as shown in Lemma~\ref{lemma:DT_Origin} and Lemma~\ref{lemma:DT_Origin_Oscill}.

Finally, notice that the following theorem holds also for structurally balanced networks. However, in that case the condition $\pi_1<\pid$ would always be satisfied (see Proposition~\ref{prop:fpi1pi1d}) meaning that the formulation of theorem could be simplified.
\begin{theorem}\label{thm:DT_summary}
Consider the system~\eqref{eqn:System_DT} where each nonlinear function $\psi_i(\cdot)$, $i=1,\mydots,n$, satisfies the properties \eqref{assumption:1psiOdd}$\div$\eqref{assumption:4Sigmoidal}. 
Assume that the signed graph $\G$ is structurally unbalanced with normalized signed Laplacian $\Lnorm$.
Let $J_\pi$, $\Lpi$, $\pi_1$ and $\pid$ be as in \eqref{eqn:Jpi}, \eqref{eqn:Lpi}, \eqref{eqn:pi1} and \eqref{eqn:pi1d}, respectively. Assume that $1-\varepsilon \max_i \delta_i \ge 0 $.
\begin{enumerate}[label=(\roman*)]
	\item If $\pi <\min\{\pi_1,\pid\}$ then the origin is the unique equilibrium point of the system~\eqref{eqn:System_DT} and it is globally asymptotically stable.
		
	\item If $\pi_1<\pid$ and the biggest eigenvalue of $J_{\pi_1}$, $\lambda_n(J_{\pi_1})=+1$, is simple, when $\pi = \pi_1$ the system~\eqref{eqn:System_DT} undergoes a pitchfork bifurcation;
		
	If $\pi_1>\pid$ and the smallest eigenvalue of $J_{\pid}$, $\lambda_1(J_{\pid})=-1$, is simple, when $\pi = \pid$ the system~\eqref{eqn:System_DT} undergoes a period-doubling bifurcation.
\end{enumerate}
\end{theorem}
\noindent%
Proof in Appendix~\ref{proof:DT_summary}.

Notice that, compared with Theorems~\ref{thm:CT_SB_summary} and \ref{thm:CT_SUB_summary}, Theorem~\ref{thm:DT_summary} considers only the first bifurcation the system~\eqref{eqn:System_DT} undergoes at the origin, i.e., it does not consider for instance secondary bifurcations at the origin happening for values of $\pi$ such that $\lambda_{n-1}(J_\pi)=+1$ or $\lambda_{2}(J_\pi)=-1$.
\begin{corollary}\label{thm:DT_BoundedEq}	
Consider the system~\eqref{eqn:System_DT} where each nonlinear function $\psi_i(\cdot)$, $i=1,\mydots,n$, satisfies the properties \eqref{assumption:1psiOdd}$\div$\eqref{assumption:4Sigmoidal}. 
Let $\frustration$ be the frustration of the signed network $\G$ defined in \eqref{eqn:frustration}.
\begin{enumerate}[label=(\roman*)]
	\item If $x^*$ is an equilibrium point of \eqref{eqn:System_DT}, then $\norm{x^*}_1\le\pi(n-2\frustration)$.
			
	\item If $\Delta = \delta I$ with $\delta \varepsilon<1$, the trajectories of \eqref{eqn:System_DT} asymptotically converge to the set $\{x\in \R^n: \norm{x}_1\le \pi n\}$.
	\end{enumerate}
\end{corollary}
The proof is omitted since (i) follows from the observation that the discrete- and continuous-time models share the same equilibrium points, therefore the upper bound on the $1$-norm of the equilibria found in Theorem~\ref{thm:CT_BoundedEq}(i) still holds, and (ii) follows from the fact that the nonlinearities are saturated (and can be shown for instance using the Lyapunov function $V(x_k)=\norm{x_k}_1-\pi n$ for all $\norm{x_k}_1>\pi n$ and $V(x_k)=0$ otherwise).

\section{Discussion and interpretation of the results}
\label{sec:Results}
\begin{figure}[t]\centering
	\includegraphics[width=.3\textwidth]{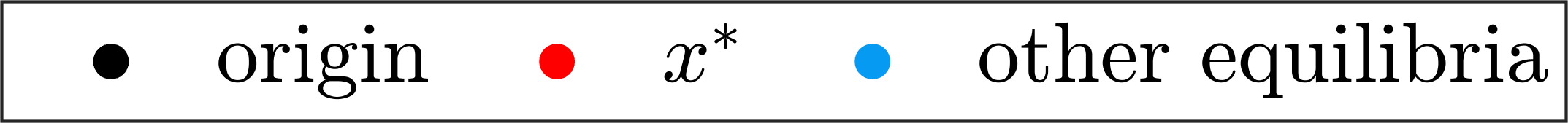}\\\vspace{-.1cm}
	\subfloat[]{\includegraphics[width=.23\textwidth]{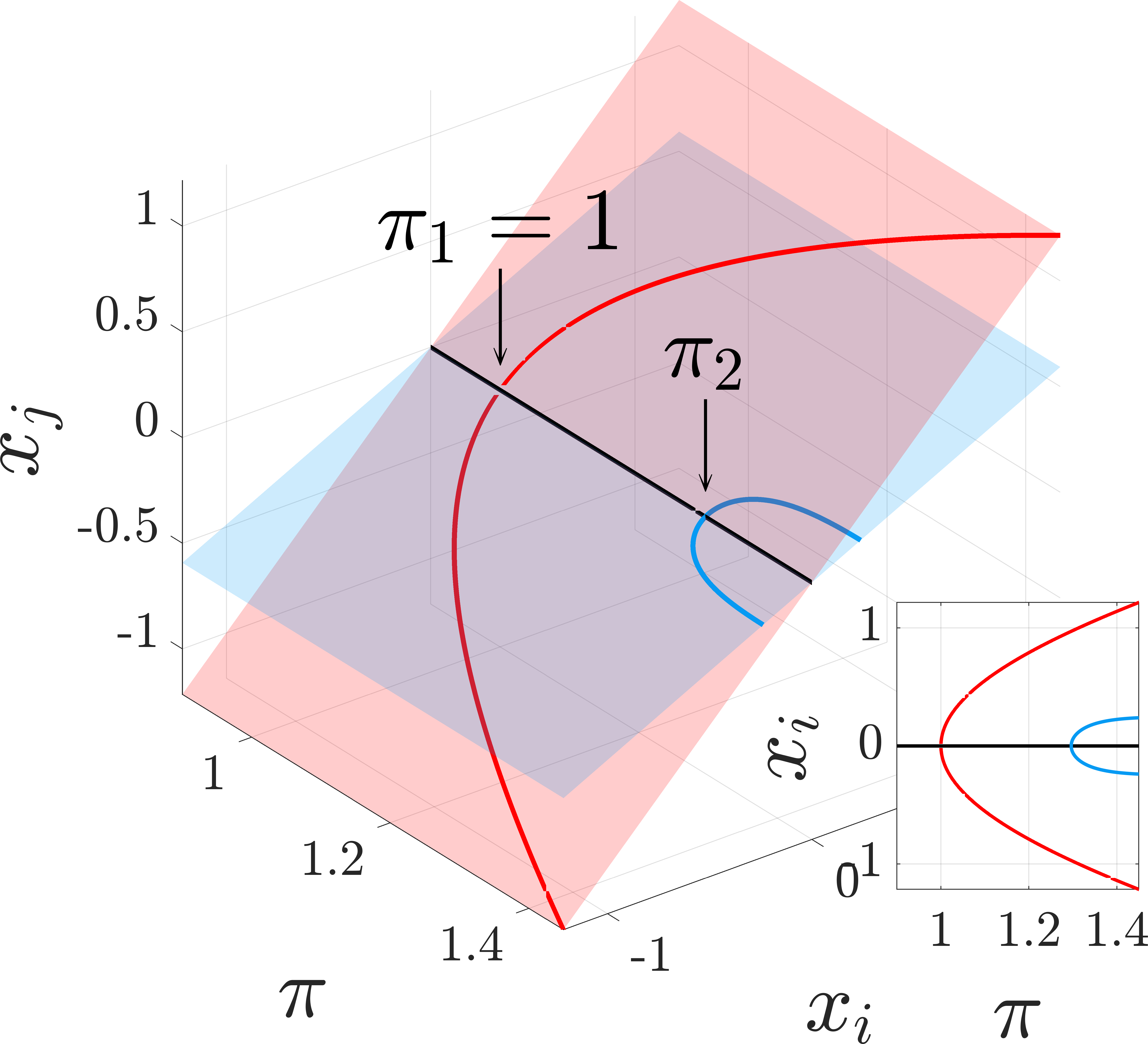}\label{fig:EX_Bifurcation_SB}}\,
	\subfloat[]{\includegraphics[width=.23\textwidth]{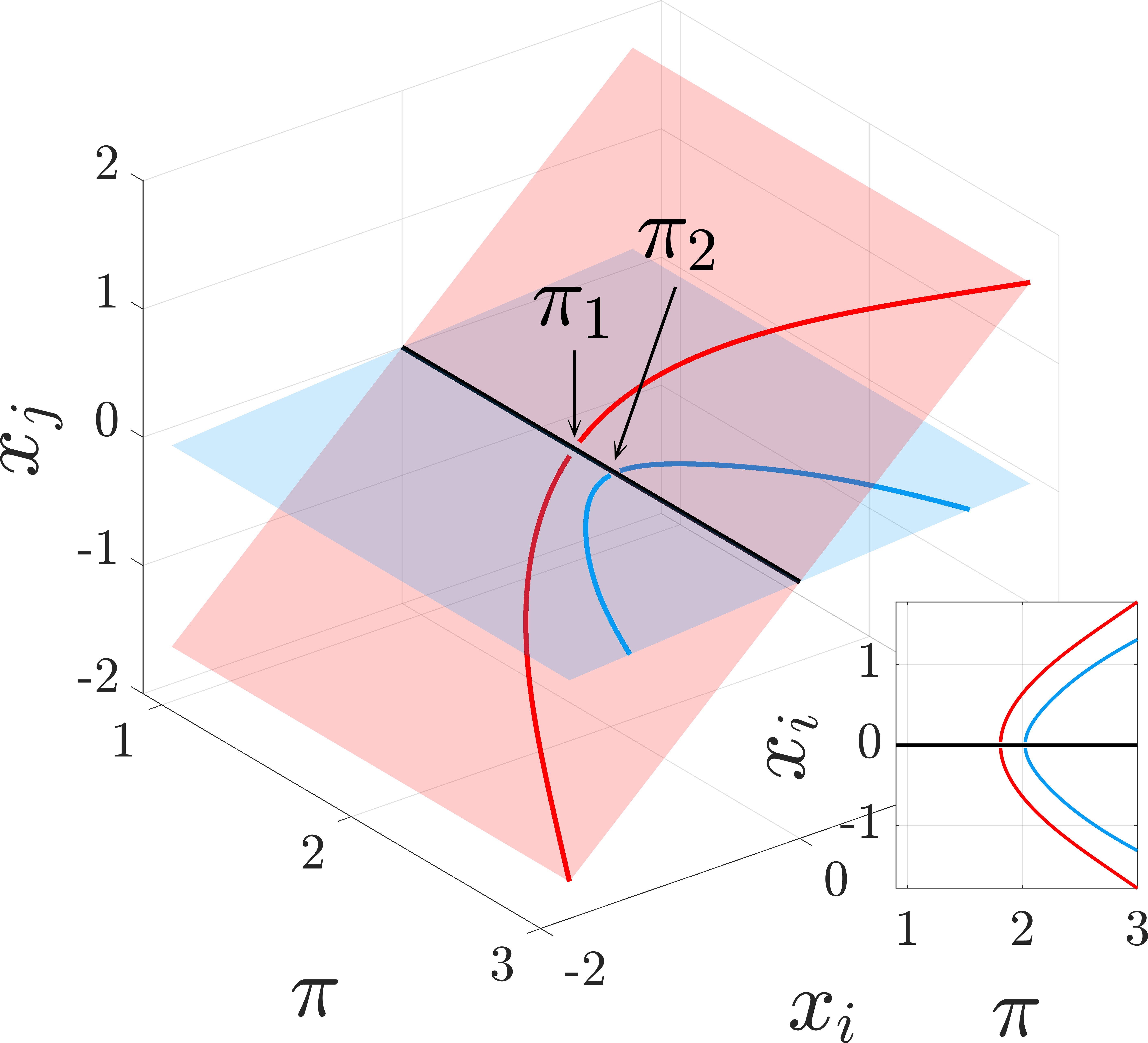}\label{fig:EX_Bifurcation_SUB}}
	\caption{The system~\eqref{eqn:System_CT} undergoes two pitchfork bifurcations, respectively for $\pi = \pi_1$ and $\pi = \pi_2$. The bifurcation diagram for two components $x_i$ and $x_j$ is here shown for two different signed networks. (a): Structurally balanced network (monotone system). (b): Structurally unbalanced network.}
	\label{fig:EX_CT_Bifurcation}
\end{figure}

Comparing Theorem~\ref{thm:CT_SB_summary} with Theorem~\ref{thm:CT_SUB_summary}, the general behavior of the continuous-time system~\eqref{eqn:System_CT} (illustrated in Figure~\ref{fig:EX_CT_Bifurcation}) does not change when, instead of a structurally balanced network, we assume that the network is structurally unbalanced. 
However, while in the structurally balanced case (see Figure~\ref{fig:EX_Bifurcation_SB}) the first threshold value for the social effort parameter $\pi$ is constant ($\pi_1=1$), in the structurally unbalanced case (see Figure~\ref{fig:EX_Bifurcation_SUB}) this value is strictly greater than $1$ and increases with the smallest eigenvalue of the normalized signed Laplacian, $\lambda_1(\Lnorm)$. 
In a recent work \cite{FontanAltafini2018CDC} we have shown that $\lambda_1(\Lnorm)$ approximates well the frustration of a signed network $\G$ (measured by $\frustration$ introduced in equation~\eqref{eqn:frustration}), while the intuition is that $\lambda_2(\Lnorm)$ is independent from the frustration $\frustration$.
As a consequence, a higher frustration $\frustration$ (reflecting a situation where the system~\eqref{eqn:System_CT} is ``far'' from being monotone) implies (i) a higher value of $\pi_1$ 
and (ii) the shrinkage of the interval $(\pi_1,\pi_2)$ for which only two alternative equilibria are admitted. These conclusions are illustrated in Example~\ref{example:Norm}.

\begin{figure}[t]\centering
	\includegraphics[width=.41\textwidth]{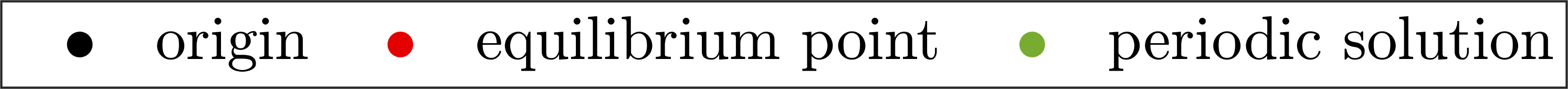}\\\vspace{-.1cm}
	\subfloat[]{\includegraphics[width=.15\textwidth]{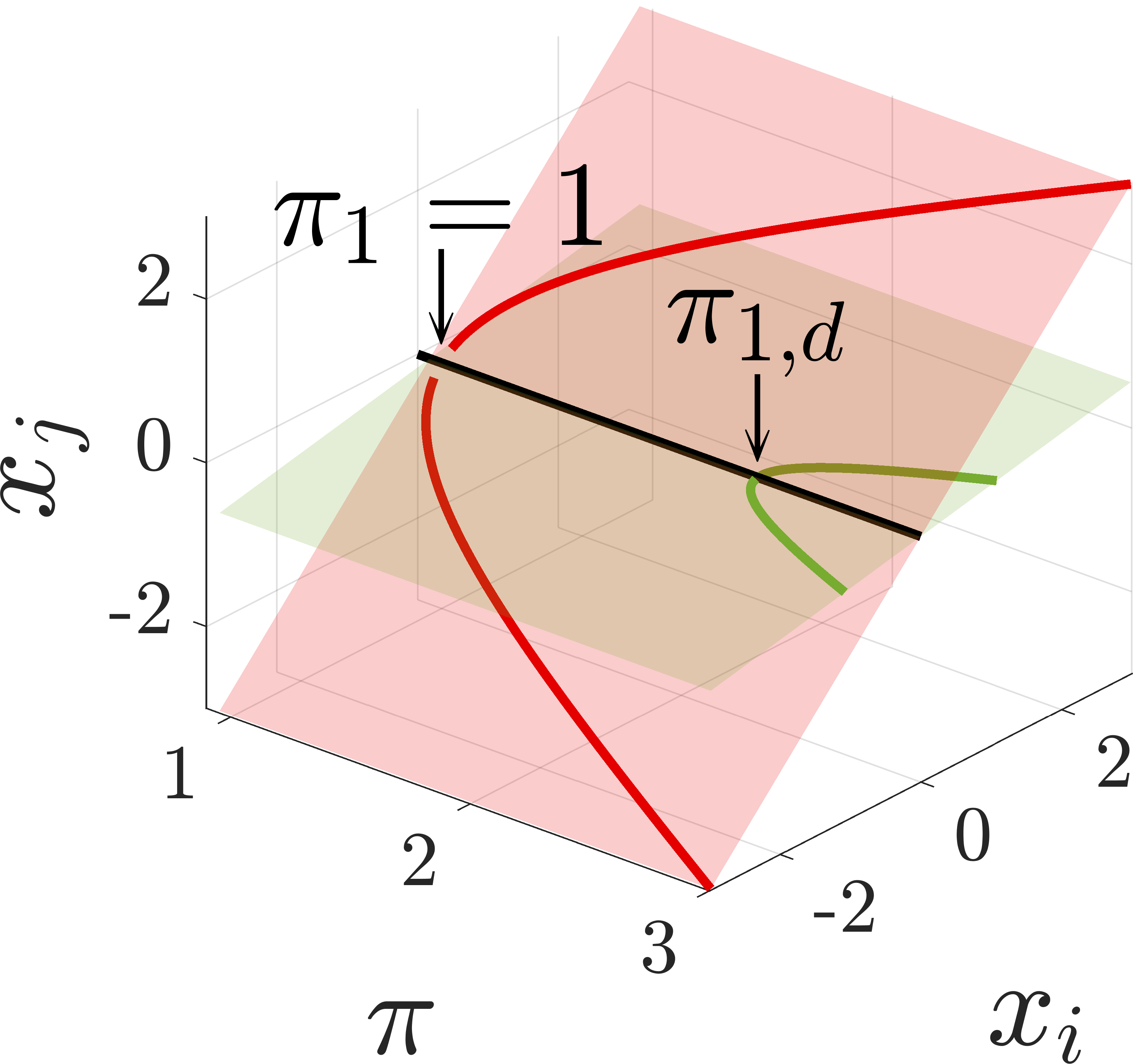}\label{fig:EX_DT_SB}}\,
	\subfloat[]{\includegraphics[width=.15\textwidth]{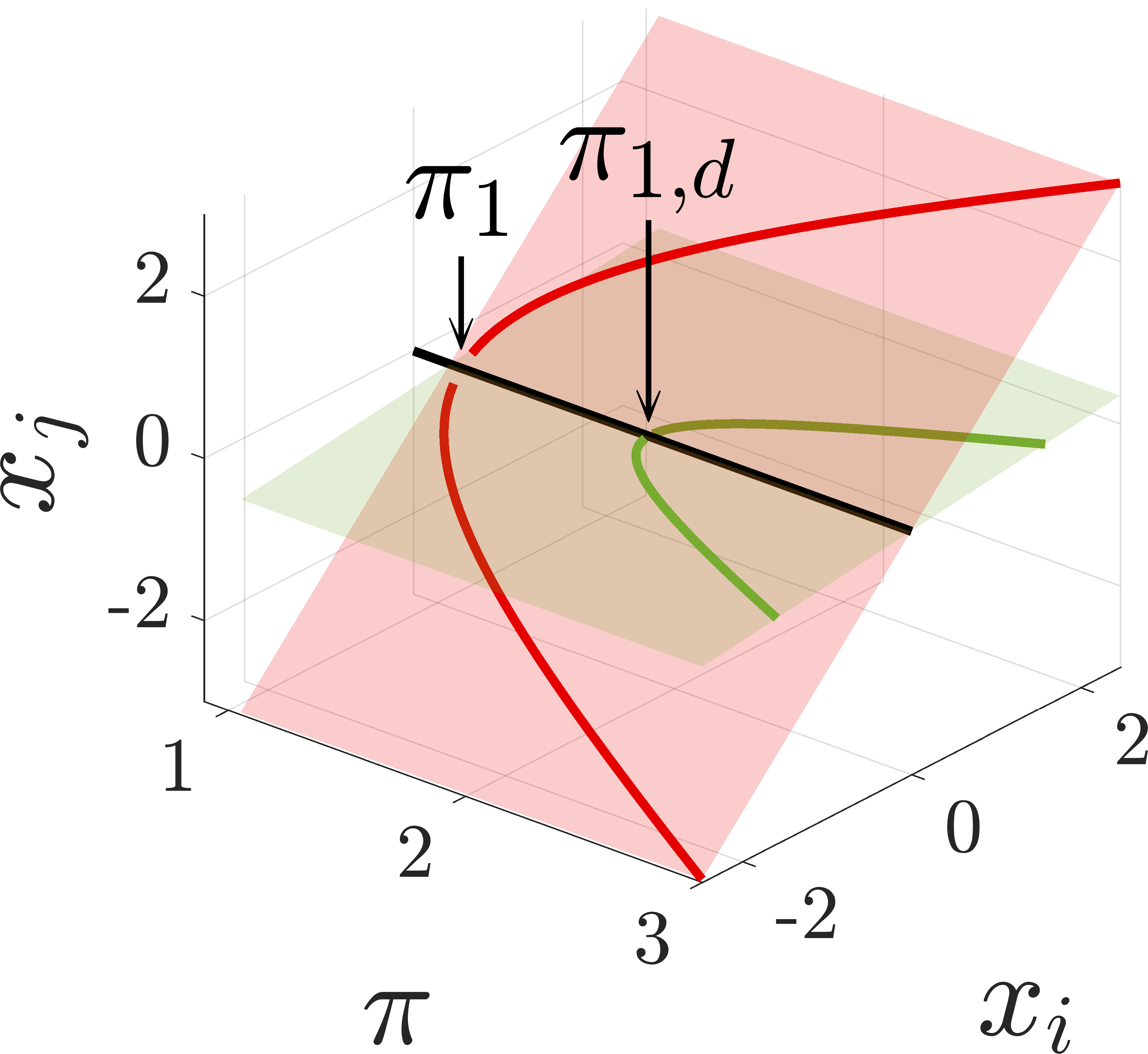}\label{fig:EX_DT_SUB1}}\,
	\subfloat[]{\includegraphics[width=.15\textwidth]{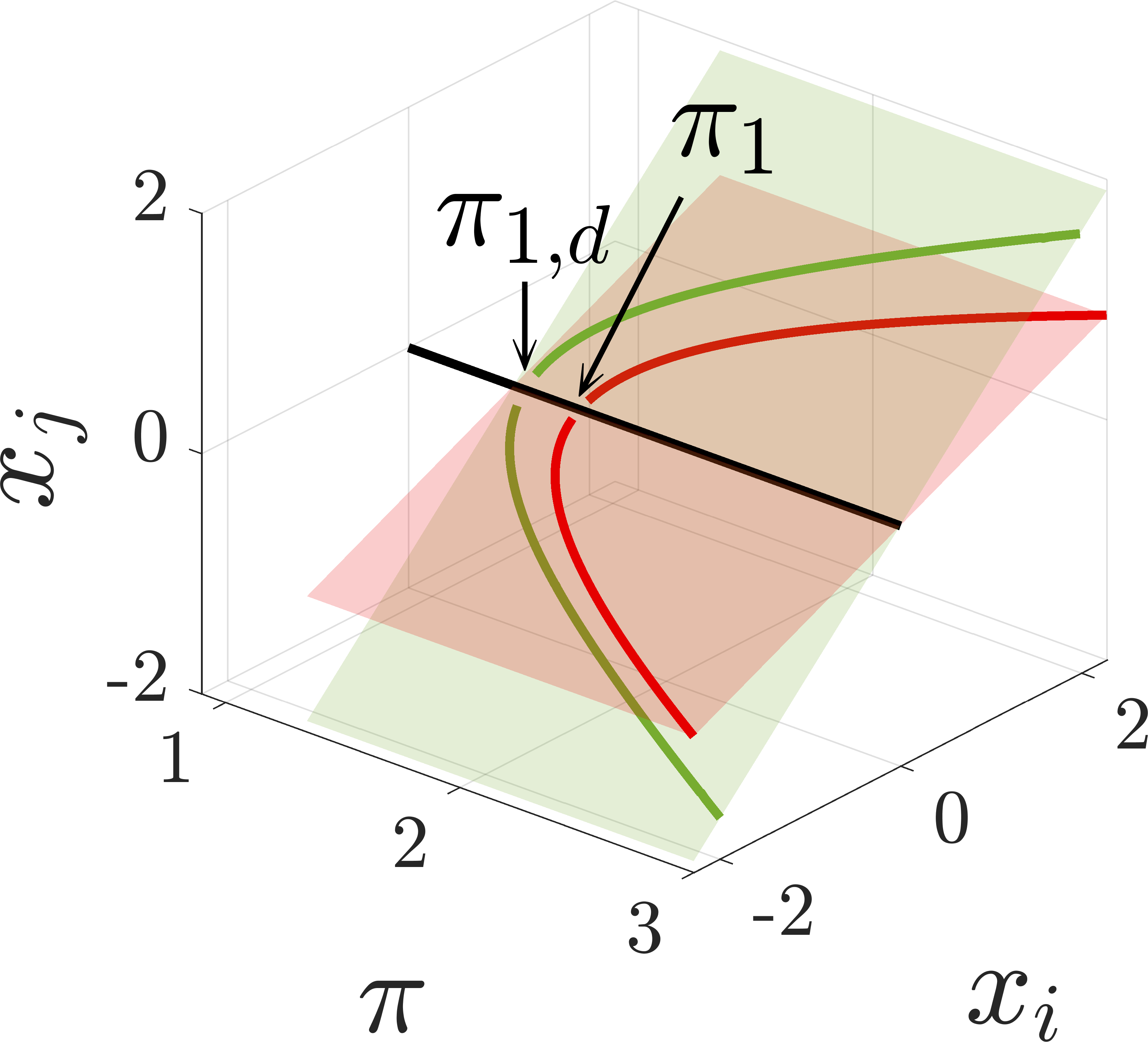}\label{fig:EX_DT_SUB2}}
	\caption{The system~\eqref{eqn:System_DT} undergoes a pitchfork bifurcations and a period-doubling bifurcation, respectively for $\pi = \pi_1$ and $\pi = \pid$. The bifurcation diagram for two components $x_i$ and $x_j$ is here shown for three different signed networks. (a): Structurally balanced network. (b): Structurally unbalanced network with $\pi_1<\pid$. (c): Structurally unbalanced network with $\pi_1>\pid$.}
	\label{fig:EX_DT_Bifurcation}
\end{figure}

In the context of social networks the decision-making process~\eqref{eqn:System_A} can be summarized as follows:
\begin{itemize}
	\item $\pi < \pi_1$: No decision will be reached if the social effort among the agents is small.
	
	\item $\pi\in (\pi_1,\pi_2)$: The ``right'' level of commitment among the agents leads to two possible (alternative) decisions. If the signed social network is not structurally balanced, a higher frustration implies that a higher effort will be required from the agents in order to achieve this decision. 
	
	\item $\pi > \pi_2$: An overcommitment between the agents (high value of social effort) leads to a situation where several alternative decisions are possible.
\end{itemize}
When we instead compare Theorems~\ref{thm:CT_SB_summary} and \ref{thm:CT_SUB_summary} with Theorem~\ref{thm:DT_summary}, we observe that the discrete-time system~\eqref{eqn:System_DT} exhibits a ``richer'' behavior, in that it admits (stable) periodic solutions, as illustrated in Figure~\ref{fig:EX_DT_Bifurcation} (which, for the sake of simplicity, does not consider secondary pitchfork or period-doubling bifurcations at the origin).
This is related to the presence of a new threshold value for the parameter $\pi$, denoted $\pid$:
understanding where $\pid$ lies compared with $\pi_1$ plays a key role when investigating the behavior of the system~\eqref{eqn:System_DT} over a signed network.
In particular, Proposition~\ref{prop:fpi1pi1d} suggests that the condition $\pi_1>\pid$ cannot hold unless a signed network is structurally unbalanced and has high frustration (i.e., $\lambda_1(\Lnorm)\gg 0$).

This implies first that, if we consider networks that are structurally balanced (for which $\pi_1<\pid$ always holds) or that are structurally unbalanced for which $\pi_1<\pid$ (typically, with low frustration) the general behavior of the discrete-time system~\eqref{eqn:System_DT} resembles that of its continuous-time counterpart, see Fig.~\ref{fig:EX_DT_SB} and Fig.~\ref{fig:EX_DT_SUB1}: the crossing of a (pitchfork) bifurcation yields two (alternative) nontrivial equilibrium points representing two possible (alternative) decisions.
Hence, the general idea that the higher is the frustration of the network the higher is the social effort needed to converge to a nontrivial equilibrium point still holds.
Instead, if we consider networks that are structurally unbalanced for which $\pi_1>\pid$ (typically, with high frustration), see Fig.~\ref{fig:EX_DT_SUB2}, then there exists an interval of values for the social effort parameter, $(\pid,\pi_1)$, for which the collective decision-making process still ends in a deadlock situation where the opinions of the agents do not settle but keep fluctuating: only by further increasing the commitment among the agents the process can be settled and the community can reach a decision.
In conclusion, in the discrete-time model the presence of high frustration in the graph leads to agents who will never cease to change opinion (a somewhat artificial behavior).
Other recent works in the literature that propose models characterized by fluctuations of opinions of the agents are for instance \cite{Acemoglu2013,Cisneros-VelardeChanBullo2019}.

\section{Numerical Examples}
\label{Examples}
In this section we first illustrate the bound~\eqref{eqn:fpi1pi2} in Proposition~\ref{prop:fpi1pi2} (Example~\ref{ex:fpi1pi2}). As a byproduct we observe (numerically) that the bound is tight if the smallest and largest eigenvalues of $\Lnorm$ satisfy the condition $\lambda_1(\Lnorm)<2-\lambda_n(\Lnorm)$ (typically, if the network does not have high frustration).
Then, we show the behavior of the system~\eqref{eqn:System_CT} over signed (structurally unbalanced) networks with increasing frustration (Example~\ref{example:Norm} and Example~\ref{example:Norm_Frustration}). In Example~\ref{example:Multistability} we show that when the social effort parameter $\pi$ crosses the second threshold $\pi_2$ the system admits multiple equilibria which are stable (i.e., several decision states for the community are possible).
Example~\ref{example:MultipleEigenvalues} is used to illustrate a case which has not been treated by our analysis. Indeed, in Example~\ref{example:MultipleEigenvalues} we illustrate the behavior of the system~\eqref{eqn:System_CT} in presence of symmetries implying an algebraic multiplicity of $\lambda_1(\Lnorm)$ higher than $1$. The case where the smallest eigenvalue of $\Lnorm$ is not simple is in fact not covered by Theorem~\ref{thm:CT_SUB_summary}. However the intuition, supported by the reading of \cite{GolubitskyStewart2002,GolubitskyStewartSchaeffer1988}, is that when $\pi>\pi_1$ the system admits multiple (more than three) equilibria.
Finally, in Example~\ref{example:DT} we illustrate the behavior of the discrete-time system~\eqref{eqn:System_DT} and compare it with that of the continuous-time system~\eqref{eqn:System_CT}.

If not specified otherwise, we assume that each nonlinear function $\psi_i(\cdot)$ ($i=1,\mydots,n$) is given by the hyperbolic tangent $\psi(\varepsilon)= \tanh(\varepsilon)$. 
Moreover, to compute numerically the frustration $\frustration$ of a signed network $\G$ we use the algorithm proposed in \cite{IacRamSorAlt2010}.
\begin{figure}[t]\centering
	\subfloat[]{\includegraphics[width=0.23\textwidth]{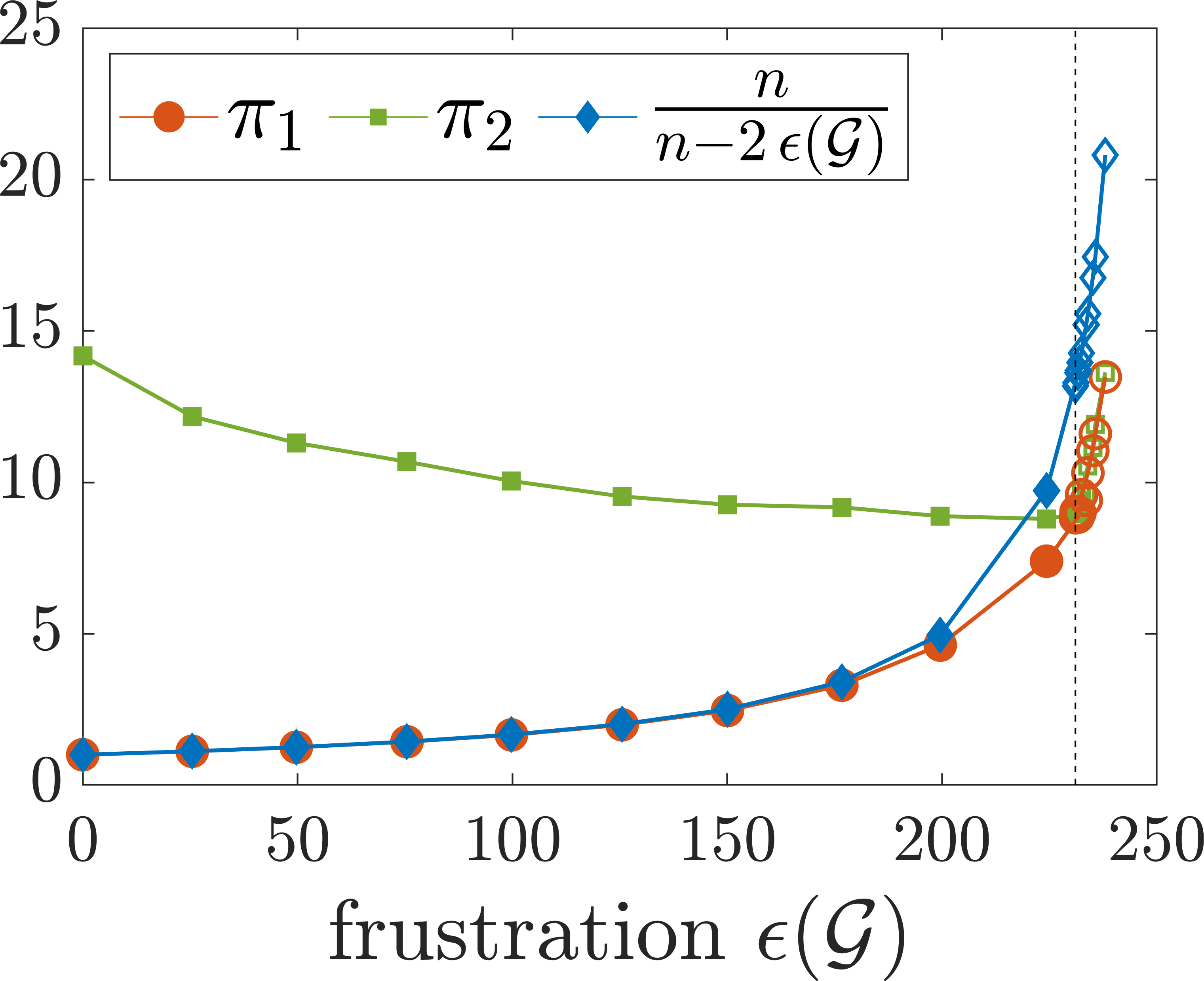}\label{fig:fpi1pi2_symm}}\,
	\subfloat[]{\includegraphics[width=0.23\textwidth]{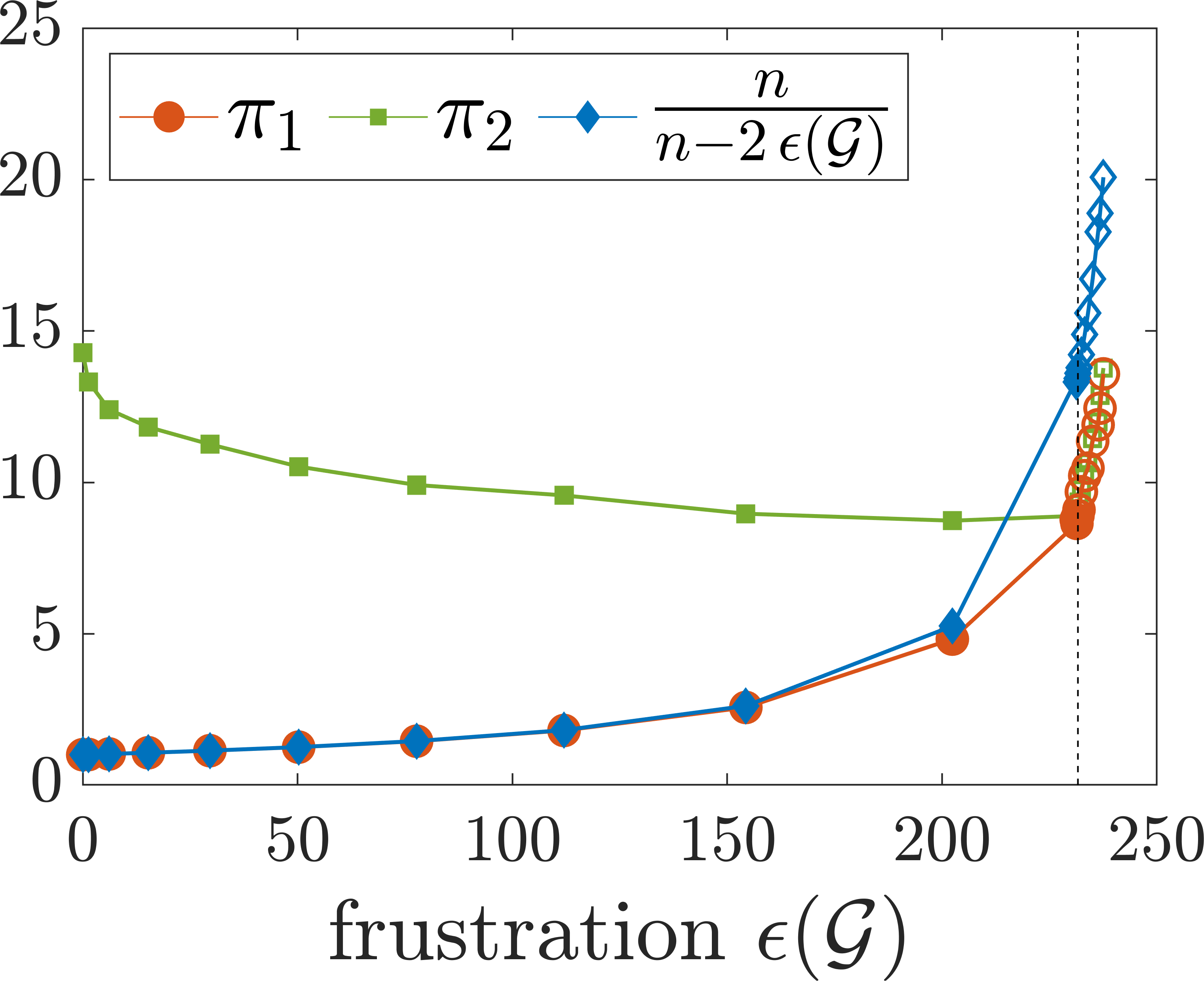}\label{fig:fpi1pi2_nosymm}}
	\caption{Example~\ref{ex:fpi1pi2}. Plot of $\pi_1,\pi_2$ and $\frac{n}{n-2\frustration}$, for two sequences of signed networks with increasing frustration $\frustration$. (a): Sequence 1: for each network, the normalized signed Laplacian $\Lnorm$ is symmetric. (b): Sequence 2: for each network, the matrix $\Lnorm$ is not symmetric. A full (resp., empty) symbol means that $\lambda_1(\Lnorm)\!<\!2\!-\!\lambda_n(\Lnorm)$ (resp., $\lambda_1(\Lnorm)>2-\lambda_n(\Lnorm)$); for clarity, a dashed line shows the maximum value of frustration above which the condition $\lambda_1(\Lnorm)<2-\lambda_n(\Lnorm)$ does not hold.}
	\label{fig:fpi1pi2}
\end{figure}
\begin{example}\label{ex:fpi1pi2}
This example wants to illustrate the bound~\eqref{eqn:fpi1pi2} in Proposition~\ref{prop:fpi1pi2} and show that it holds also for graphs whose normalized signed Laplacian is not symmetric.
In Fig.~\ref{fig:fpi1pi2} we consider two sequences of signed networks $\G$ with $n\!=\!500$ agents (in which the edge weights are drawn from a uniform distribution and $p=0.8$ is the edge probability) and with increasing frustration $\frustration$.
In the first sequence (see Fig.~\ref{fig:fpi1pi2_symm}), each adjacency matrix $A$ of the network is rescaled so that $\abs{A}\1= \delta\1$, which implies that the normalized signed Laplacian $\Lnorm$ is symmetric.
Instead, in the second sequence (see Fig.~\ref{fig:fpi1pi2_nosymm}), each matrix $\Lnorm$ is not symmetric (but is symmetrizable).
As Fig.~\ref{fig:fpi1pi2} illustrates, the bound~\eqref{eqn:fpi1pi2} holds for both sequences; moreover, when the frustration is small (numerically, when the condition $\lambda_1(\Lnorm)<2-\lambda_n(\Lnorm)$ is satisfied) the upper bound $\frac{n}{n-2\frustration}$ for $\pi_1$ is tight (this is not surprising since we know that $\lambda_1(\Lnorm)$ approximates well the frustration).
\end{example}

\begin{figure}[t]\centering
	\includegraphics[width=.28\textwidth]{EX_Bifurcation_legend}\\\vspace{-.1cm}
	\subfloat[]{\includegraphics[height=.15\textwidth]{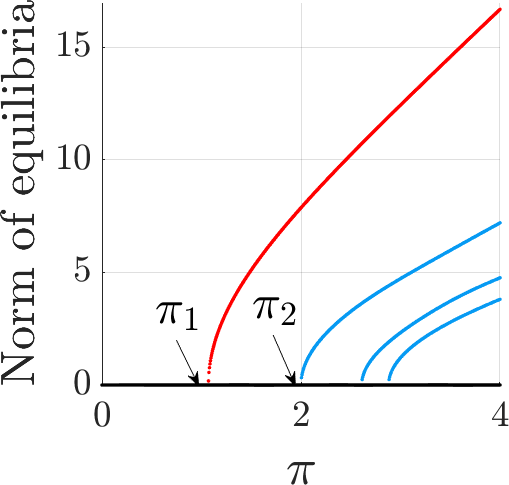}\label{fig:EX_Norm_1}}\,
	\subfloat[]{\includegraphics[height=.15\textwidth]{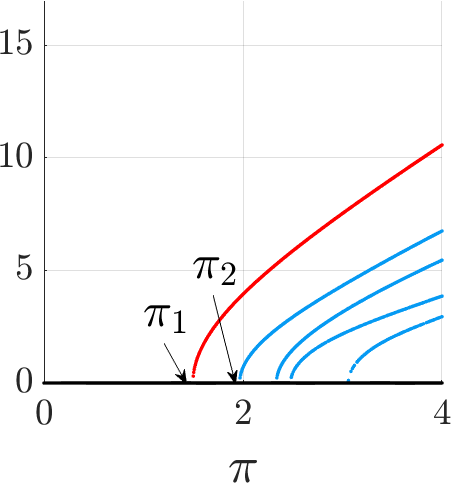}\label{fig:EX_Norm_2}}\,
	\subfloat[]{\includegraphics[height=.15\textwidth]{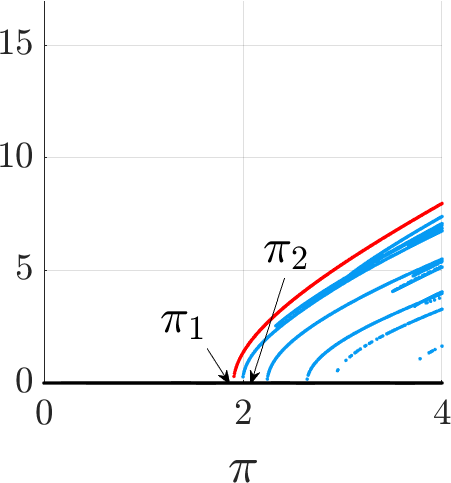}\label{fig:EX_Norm_3}}
	\caption{Example~\ref{example:Norm}. Norm of the equilibrium points of the system~\eqref{eqn:System_CT} as a function of $\pi$. The networks $\G$ we use in this example are structurally unbalanced, with increasing frustration $\frustration$, see Table~\ref{table:EX_Norm}. 
		(a): $ \frustration = 0.677 $. 
		(b): $ \frustration = 4.285 $. 
		(c): $ \frustration = 5.536 $. 
	}
	\label{fig:EX_Norm}
\end{figure}

\begin{example}
	\label{example:Norm}
	We consider three signed networks with $n=20$ agents in which the edge weights are drawn from a uniform distribution and $p=0.5$ is the edge probability. These networks are chosen to be structurally unbalanced and with increasing frustration.
	In Fig.~\ref{fig:EX_Norm} the euclidean norm of the equilibria of the system~\eqref{eqn:System_CT} for values of $\pi$ in $\{0.005,0.010,\mydots,4\}$ is depicted. As Table~\ref{table:EX_Norm} shows, the smallest eigenvalue of the normalized signed Laplacian increases with the frustration of the network while the second smallest eigenvalue remains almost constant, hence the interval for $\pi$ for which the system admits only two equilibria becomes smaller (compare Fig.~\ref{fig:EX_Norm_1} and Fig.~\ref{fig:EX_Norm_3}).
\end{example}
\begin{table}[t]\centering	
	\begin{tabular}{Z|Z Z Z Z Z}
		&\frustration & \lambda_1 & \lambda_2 & \pi_1 & \pi_2\\\toprule
		\text{(a)} & 0.677 & 0.065 & 0.500 & 1.069 & 2.000 \\
		\text{(b)} & 4.285 & 0.332 & 0.491  & 1.496 & 1.966 \\
		\text{(c)} & 5.536 & 0.475 & 0.499 & 1.905 & 1.995\\\bottomrule
	\end{tabular}
	\caption{Example~\ref{example:Norm}. Values of frustration, first two eigenvalues of the normalized signed Laplacian and bifurcation points for the three cases, (a), (b) and (c), depicted in Figure~\ref{fig:EX_Norm}.}
	\label{table:EX_Norm}
\end{table}

\begin{figure}[!ht]\centering
	\includegraphics[width=0.45\textwidth]{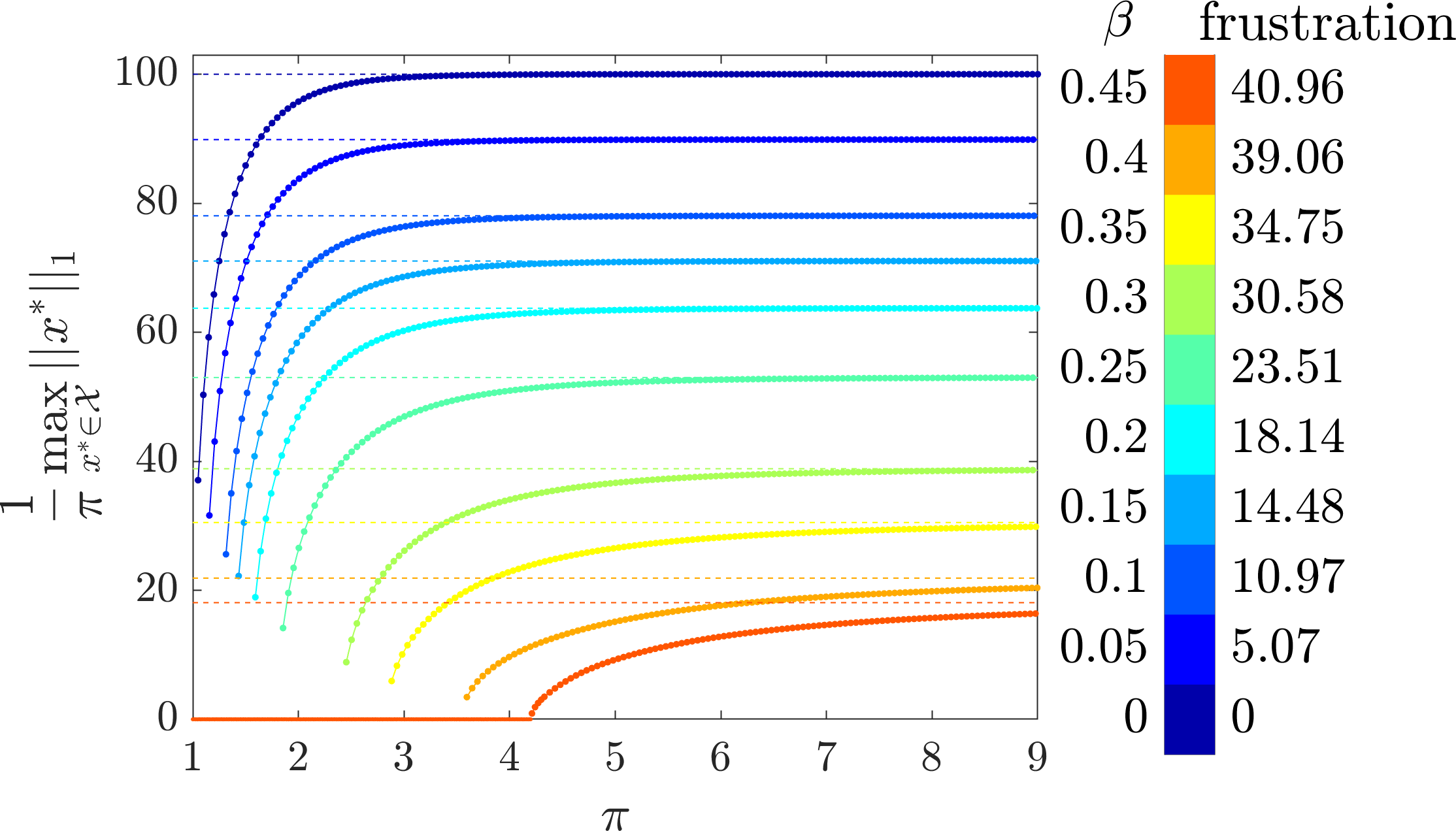}
	\caption{Example ~\ref{example:Norm_Frustration}. Plot of the maximum 1-norm of $x^*$, where $x^\ast$ is an equilibrium point of the system~\eqref{eqn:System_CT}, for a sequence of signed networks with increasing frustration. The values of $n-2 \epsilon(\G_\beta)$ are shown as dotted lines.}
	\label{fig:norm1_comparison_frustration}
\end{figure}

\begin{example}
	\label{example:Norm_Frustration}
	Consider a network $\G$ with $n=100$ agents in which the edge weights are drawn from a uniform distribution and $p=0.8$ is the edge probability.
	Let $A=[a_{ij}]$ be its weighted nonnegative adjacency matrix.
	Consider now a sequence of signed networks $\G_\beta$ with weighted adjacency matrices $A_\beta=[{a_\beta}_{ij}]$ constructed such that $\abs{A_\beta}=A$ and their signature is dependent on a parameter $\beta \in \{0,0.05,0.1,\mydots,1\}$:
	if $a_{ij}\ne 0$ then ${a_\beta}_{ij}\ne 0$ and $P[{a_\beta}_{ij}<0] = \beta$.
	When $\beta = 1$, $A_1 = -\abs{A}$.
	As $\beta$ increases also the frustration of the networks increases. 
	
	For each network, we numerically compute the equilibria $x^\ast$ of the system~\eqref{eqn:System_CT} for values of $\pi$ in $\{1,1.05,\mydots,9\}$ and their 1-norm $\norm{x^\ast}_1$: let $\mathcal{X}= \{x^\ast\in \R^n$: $x^\ast$ is an equilibrium point of the system~\eqref{eqn:System_CT}$\}$ be the set of equilibria.	
	In Fig.~\ref{fig:norm1_comparison_frustration}, for each network of the sequence we plot $\displaystyle\frac{1}{\pi} \max_{x^\ast\in \mathcal{X}} \norm{x^\ast}_1$ (the maximum 1-norm of the equilibrium points divided by $\pi$) for each value of $\pi$; the colormap illustrates the sequence of signed networks $\G_\beta$ with increasing frustration.
	As Theorem~\ref{thm:CT_BoundedEq} states, the maximum 1-norm of the equilibria is upper bounded by $\pi(n-2\epsilon(\G_\beta))$, where $\epsilon(\G_\beta)$ indicates the frustration of $\G_\beta$: as the frustration increases, the bound decreases.
\end{example}
\begin{figure}[!ht]\centering
	\includegraphics[width=.35\textwidth]{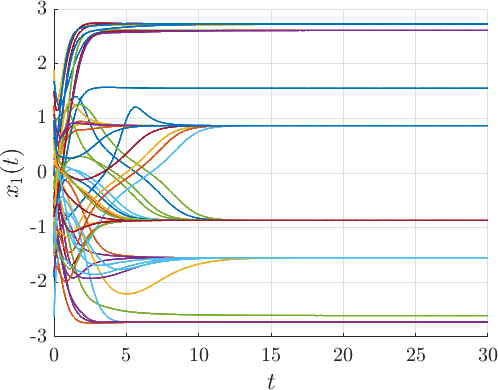}
	\caption{Example~\ref{example:Multistability}. Evolution of state variable $x_1(t)$ for $50$ random initial conditions and $\pi = 4$. The signed network considered in this example corresponds to the one used to obtain Fig.~\ref{fig:EX_Norm_3} ($\pi_2 = 1.995 $).}
	\label{fig:EX_Multistability}
\end{figure}
\begin{example}
	\label{example:Multistability}
	When $\pi>\pi_2$, Theorem~\ref{thm:CT_SUB_summary}(iii) proves that the system~\eqref{eqn:System_CT} admits multiple equilibrium points. Through numerical simulations it is possible to see that some of these equilibria may be stable.	
	In Fig.~\ref{fig:EX_Multistability} the multistability for the system~\eqref{eqn:System_CT} is highlighted as we depict the evolution of the first component of $x(t)$ for $50$ random initial conditions and $\pi = 4>\pi_2 = 1.995$. The same signed network as case (c) in Fig.~\ref{fig:EX_Norm} and Table~\ref{table:EX_Norm} is used.
\end{example}

\begin{figure}[t]\centering
	\subfloat[]{\includegraphics[width=.45\textwidth]{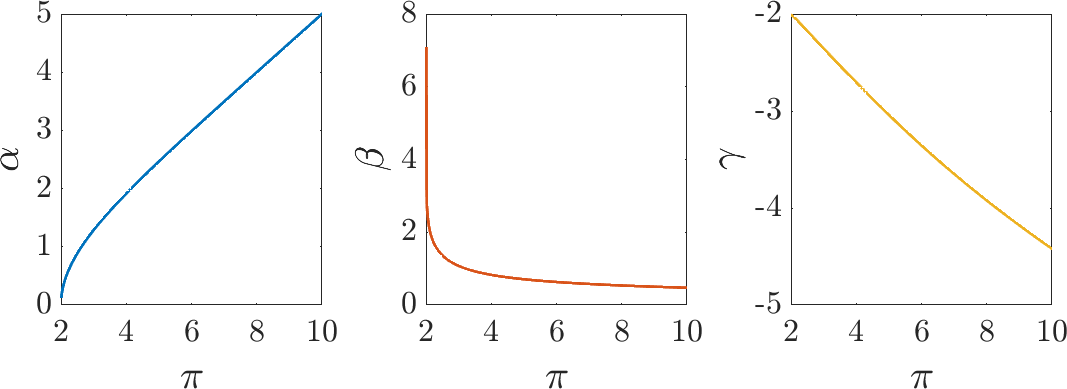}\label{fig:EX_MultipleEig_abc}}\\%
	\subfloat[]{\includegraphics[width=.37\textwidth]{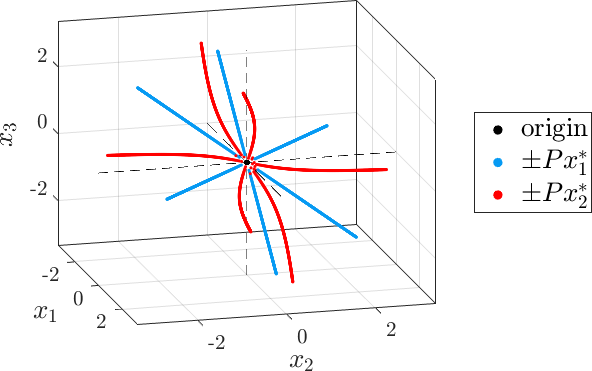}\label{fig:EX_MultipleEig_eq}}
	\caption{Example~\ref{example:MultipleEigenvalues}.
		(a): $\alpha,\beta,\gamma$ as functions of $\pi$.	
		(b): Equilibria of the system~\eqref{eqn:System_CT} as described by \eqref{eqn:EX_MultipleEig_Eq}, for $\pi = 2.001,2.002,\dots,4$: the origin (black dot), $\pm P x_1^*$ (blue branches) and $\pm P x_2^*$ (red branches).}
\end{figure}

\begin{example}
	\label{example:MultipleEigenvalues}
	Consider the system~\eqref{eqn:System_CT} where each nonlinear function $\psi_i(\cdot)$, $i=1,\mydots,n$, satisfies the properties \eqref{assumption:1psiOdd}$\div$\eqref{assumption:4Sigmoidal}.
	Moreover, assume that 
	\begin{multline}
	\psi_i(\varepsilon)=\psi_j(\varepsilon) =: \psi_u(\varepsilon),\,\forall i,j=1,\mydots,n,\;\varepsilon\in \R \\\text{(identical nonlinearities).}
	\tag{A.5}
	\label{assumption:5psiEqual}
	\end{multline}
	Notice that under these assumptions $P\psi(x) = \psi(Px)$ for all signed permutation matrices $P$.
	
	Let $n=3$ and the adjacency matrix of the network be
	\begin{equation*}
	A = 
	\begin{bmatrix} 
	0&-1&-1\\
	-1&0&-1\\
	-1&-1&0
	\end{bmatrix} = I -\1\1^T,
	\end{equation*}
	which implies that the signed graph described by $A$ is structurally unbalanced	and that the smallest eigenvalue of the normalized signed Laplacian $\Lnorm$, $\lambda_1(\Lnorm)$, is not simple (the spectrum of $\Lnorm$ is $\spectrum{\Lnorm}=\left\{ \frac{1}{2},\frac{1}{2},2\right\}$).
	This represents an interesting case for our analysis since the (algebraic and geometric) multiplicity of the smallest eigenvalue of $\Lnorm$ is $2$, hence we cannot straightforwardly apply Theorem~\ref{thm:CT_SUB_summary}(ii).
	However, in this case the equilibria of the system~\eqref{eqn:System_CT} for $\pi >\pi_1=\frac{1}{1-\lambda_1(\Lnorm)}=2$ can be computed explicitly.
	
	Under assumption \eqref{assumption:5psiEqual}, let $\pi>2$ and $\alpha(\pi),\beta(\pi)>0,\gamma(\pi)< 0$ be such that
	\begin{equation*}
		\alpha:\,\frac{\psi_u(\alpha)}{\alpha} = \frac{2}{\pi},
		\quad
		\beta,\gamma:\, \begin{cases}
		\gamma = -\pi \psi_u(\beta)\\
		\psi_u(\beta) + \frac{2}{\pi} \beta + \psi_u(\gamma) = 0,
		\end{cases}
	\end{equation*}
	see also Fig.~\ref{fig:EX_MultipleEig_abc}.
	Then $x_1^* = \alpha [1 , -1 , 0]^T$, $x_2^* =[\beta, \beta , \gamma ]^T$
	are equilibrium points of \eqref{eqn:System_CT}. Indeed
	\begin{gather*}
		\pi H \psi(x_1^*) 
		= \pi \psi_u(\alpha) H \frac{x_1^*}{\alpha}
		=\frac{\pi}{2} \frac{\psi_u(\alpha)}{\alpha} x_1^*
		= x_1^*,\\
	\pi H \psi(x_2^*) 
	= -\frac{\pi}{2} \begin{bmatrix} \psi_u(\beta)+\psi_u(\gamma) \\ \psi_u(\beta)+\psi_u(\gamma) \\2\psi_u(\beta)\end{bmatrix}
	= \begin{bmatrix}  \beta \\ \beta \\\gamma \end{bmatrix}
	= x_2^*.
	\end{gather*}
	Let $\Phi(x,\pi)= - x + \pi H \psi(x)$. 
	Under assumption~\eqref{assumption:1psiOdd}, $\Phi(x,\pi)$ is odd. Moreover, since $PHP^T=H$ for all permutation matrices $P\in \R^{3\times 3}$, it holds that
	\begin{equation*}
	P \Phi(x,\pi) = \Phi(P x,\pi)\quad \forall P\in \mathbf{S}_3,
	\end{equation*}
	that is, $\Phi(x,\pi)$ is $\mathbf{S}_3$-equivariant ($\mathbf{S}_3$ indicates the symmetric group of order $3$, i.e., the group of all permutations of a three-element set).
	Hence if $x(t)$ is a solution of \eqref{eqn:System_CT}, then $\pm P x(t)$, $P\in \mathbf{S}_3$, is also a solution of \eqref{eqn:System_CT} \cite{GolubitskyStewart2002}.
	
	To conclude, the equilibria of \eqref{eqn:System_CT} can be written as
	\begin{equation}
		\pm P x_1^*,\quad \pm P x_2^*\quad \forall\; P\in \mathbf{S}_3.
		\label{eqn:EX_MultipleEig_Eq}
	\end{equation}
	Figure~\ref{fig:EX_MultipleEig_eq} shows the equilibrium points of the system~\eqref{eqn:System_CT}, where the nonlinear function $\psi_u$ is the hyperbolic tangent $\psi_u(\varepsilon)= \tanh(\varepsilon)$,	as $\pi$ increases.
\end{example}

\begin{example}\label{example:DT}
This last example wants to illustrate the results of Theorem~\ref{thm:DT_summary} for the discrete-time system~\eqref{eqn:System_DT} and compare them with the results of Theorem~\ref{thm:CT_SUB_summary} for the continuous-time system~\eqref{eqn:System_CT}.
We consider two structurally unbalanced networks ($\G_{1}$ and $\G_{2}$) with $n=6$ agents in which the edge weights are drawn from a uniform distribution and $p=0.9$ is the edge probability.
The network $\G_{1}$ is such that $1<\pi_1=1.53<1.89=\pid$ while the network $\G_{2}$ is such that $1<\pid=1.40<1.63=\pi_1$, where $\pi_1$ and $\pid$ are defined in \eqref{eqn:pi1} and \eqref{eqn:pi1d}, respectively.
	
Figure~\ref{fig:EX_DT} plots the trajectories of the discrete-time system~\eqref{eqn:System_DT} with $\varepsilon=0.3$ (top panels) and the trajectories of the continuous-time system~\eqref{eqn:System_CT} (bottom panels) for different values of $\pi$ and the same initial condition $x(0)=[-1.51,1.81,-0.12,1.23,0.49,0.91]^T$: in Fig.~\ref{fig:EX_DT_Pitchfork} we consider the network $\G_{1}$, while in Fig.~\ref{fig:EX_DT_Periodic} the network $\G_{2}$. 
When $\pi_1<\pid$ (see Fig.~\ref{fig:EX_DT_Pitchfork}), we expect the trajectories of both the discrete- and continuous-time system to converge to the origin for all values of $\pi$ less than $\pi_1$ (see left panels) and to converge to a nontrivial equilibrium point for values of $\pi$ greater than (and in a neighborhood of) $\pi_1$ (see right panels).
When $\pi_1>\pid$ (see Fig.~\ref{fig:EX_DT_Periodic}), we expect the trajectories of both the discrete- and continuous-time system to converge to the origin for all values of $\pi$ less than $\pid$ (see left panels). However, when $\pi \in (\pid,\pi_1)$ (see middle panels), while the trajectories of continuous-time system still converge to the origin, the discrete-time system admits a periodic solution.
Finally, for both the discrete- and continuous-time system to admit a nontrivial equilibrium point $\pi$ needs to be greater than $\pi_1$ (see right panels). 
\begin{figure}\centering
	\includegraphics[width=.38\textwidth]{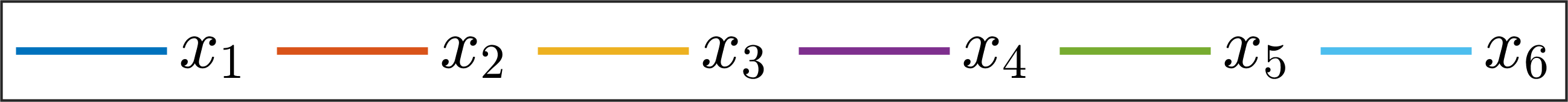}\\\vspace{-.1cm}
	\subfloat[]{\includegraphics[height=.38\textwidth]{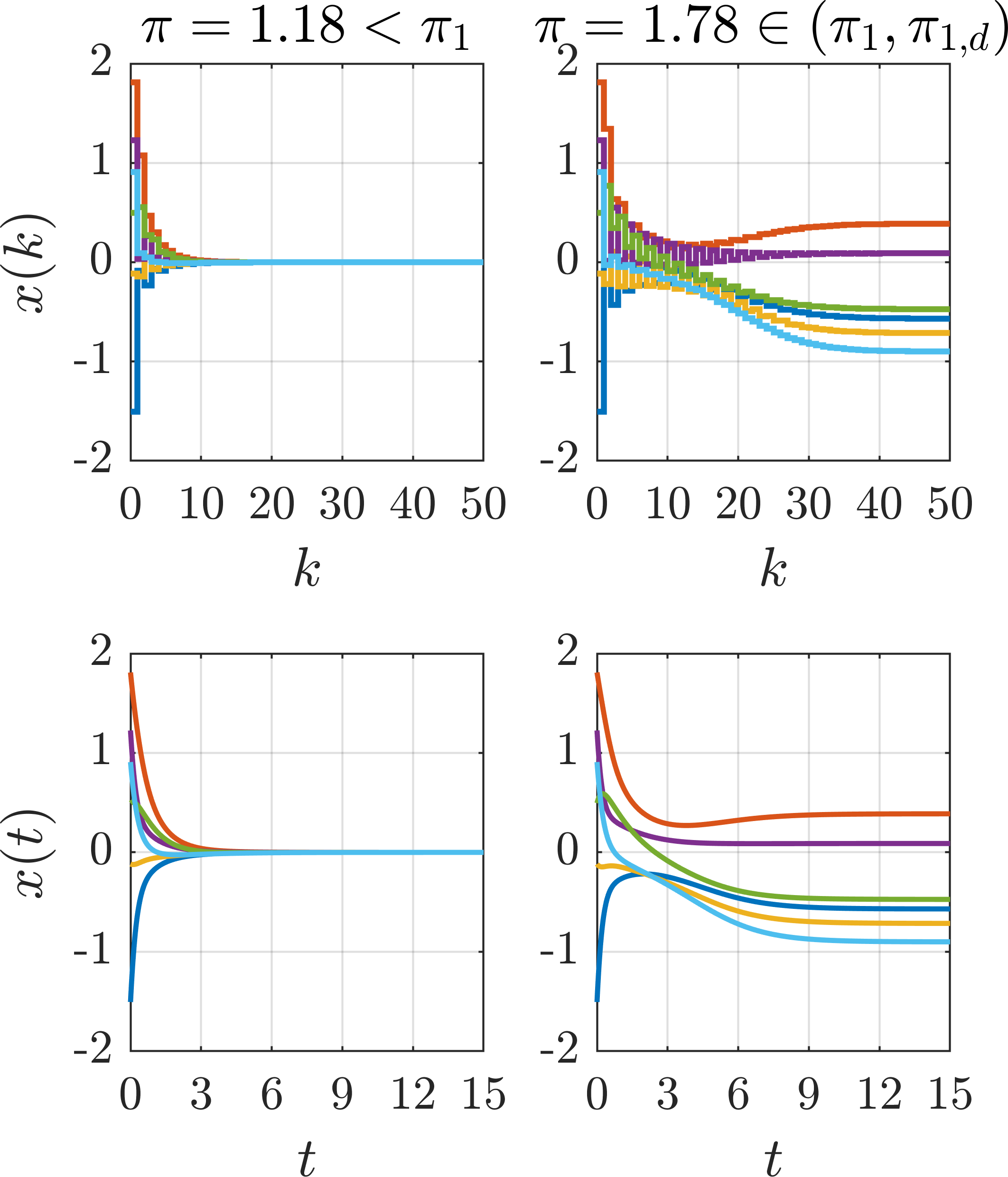}\label{fig:EX_DT_Pitchfork}}
	\\
	\subfloat[]{\includegraphics[height=.38\textwidth]{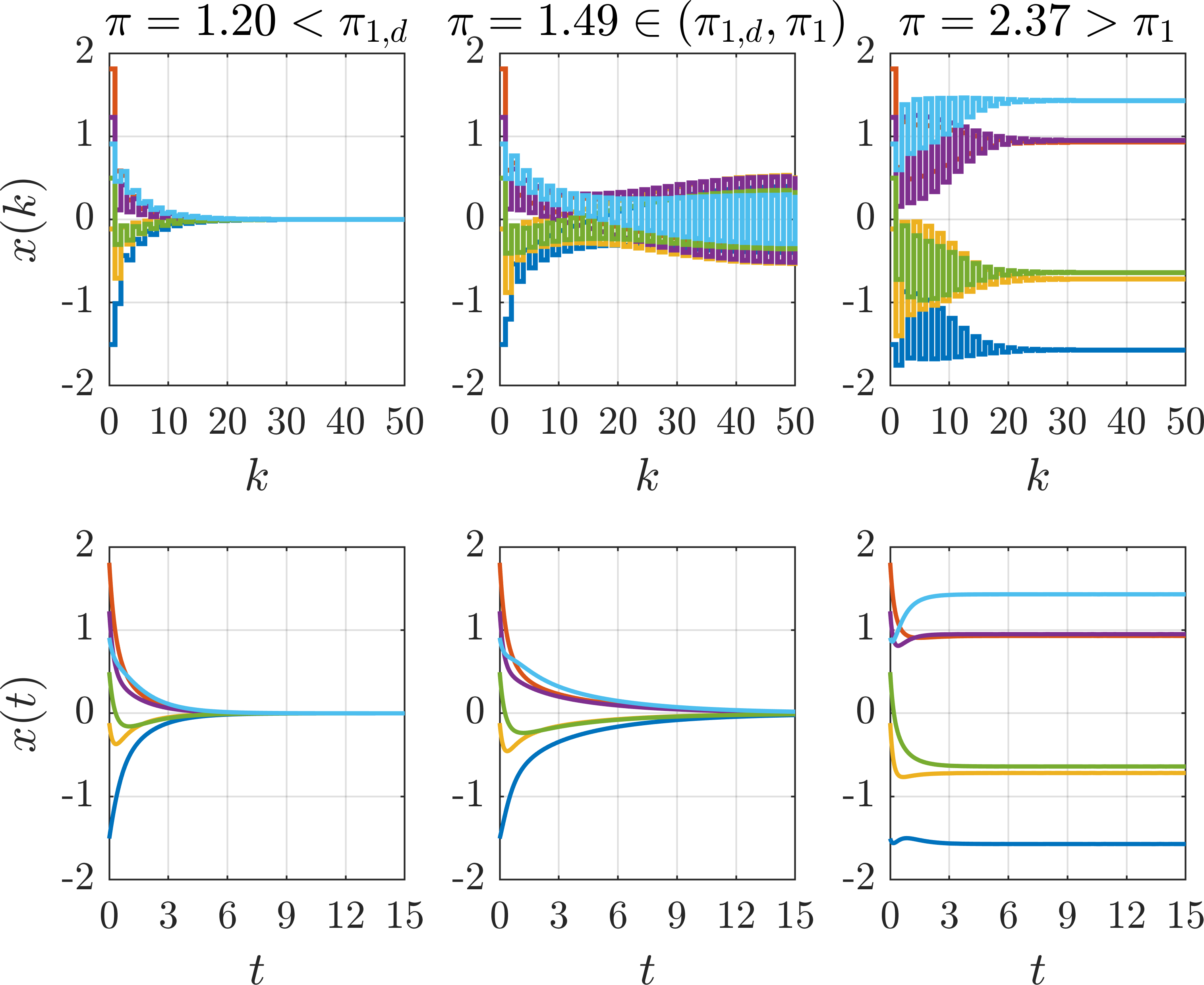}\label{fig:EX_DT_Periodic}}
	\caption{Example~\ref{example:DT}. Trajectories of the discrete-time system~\eqref{eqn:System_DT} with $\varepsilon=0.3$ (top panels) vs trajectories of the continuous-time system~\eqref{eqn:System_CT} (bottom panels) for different values of $\pi$.
	(a): Network $\G_{1}$, $\pi_1<\pid$. (b): Network $\G_{2}$, $\pi_1>\pid$.}
	\label{fig:EX_DT}
\end{figure}
\end{example}

\section{Conclusions}
\label{Conclusions}
In this work we have extended the analysis of a decision-making process in a community of agents, described by the nonlinear interconnected model introduced in \cite{GrayAl2018,FontanAltafini2018}, to the case in which the signed network representing the group of agents is not structurally balanced. 
We provided necessary and sufficient conditions for the existence (and stability) of equilibrium points of the system showing that, qualitatively, the bifurcation behavior of the system does not change when we assume that it is not monotone, i.e., that the signed social network is not structurally balanced.
What changes, however, is the threshold at which the bifurcation occurs. In particular, we have shown in the paper that this bifurcation threshold grows with the frustration of the signed network.

Given the interpretation of the bifurcation parameter as ``social effort'' of the network of agents, from a sociological point of view, this behavior is reasonable and plausible: the more ``disorder'' (i.e., frustration) a social network contains, the more difficult it is for its actors to achieve a common decision.



\appendices

\section{Technical preliminaries}
\label{sec:TechnicalPreliminaries}
In this section we introduce definitions and technical theorems and lemmas from linear algebra that will be necessary in order to prove the main results of this work.
\begin{definition}[\cite{KaszkurewiczBhaya2000,CamizStefani1996}]
	A matrix $A\in \R^{n\times n}$ is (diagonally) symmetrizable if $DA$ is symmetric for some diagonal matrix $D$ with positive diagonal entries. The matrices $DA$ and $D$ are called symmetrization and symmetrizer of $A$, respectively.
\end{definition}
\begin{theorem}[Ostrowski, 4.5.9 in \cite{HornJohnson2013}]
	\label{theorem:Ostrowski}%
	Let $A,S\in \R^{n\times n}$ with $A$ symmetric and $S$ nonsingular. Let the eigenvalues of $A$, $SAS^T$ and $SS^T$ be arranged in nondecreasing order. For each $k=1,\mydots,n$, there exists a positive real number $\theta_k$ such that $\lambda_1(SS^T)\le \theta_k \le \lambda_n(SS^T)$ and $\lambda_k(SAS^T)=\theta_k\lambda_k(A)$.
\end{theorem}
The following lemma results from Theorem~\ref{theorem:Ostrowski}.
\begin{lemma}
	\label{lemma:DiagonalEquivalence}
	Let $B\in \R^{n\times n}$ be symmetrizable, and $S = \diag{s_1,\mydots,s_n}\in \R^{n\times n}$ be a positive definite diagonal matrix. 
	Let the eigenvalues of $B$, $BS$, $SB$ and $S^{\half}BS^{\half}$ be arranged in nondecreasing order.
	Then, for all $k\in \{1,\mydots,n\}$, it holds that $\exists\,\theta_k\in \bigl[ \min_i\{s_i\},\max_i\{s_i\}\bigr]$ such that
	\begin{equation*}
	\lambda_k(BS) = \lambda_k(SB) = \lambda_k(S^{\half}BS^{\half})= \theta_k \lambda_k(B).
	\end{equation*}
\end{lemma}
\textbf{Proof.}
The matrix $B$ is symmetrizable, hence there exist a diagonal positive definite matrix $D\in \R^{n\times n}$ and a symmetric matrix $A\in \R^{n\times n}$ such that $B = D A$.  
Define the symmetric matrix
\begin{equation*}
\sym{B} : = D^{\half} A D^{\half} = D^{-\half} B D^{\half} \sim B
\end{equation*}
which, from similarity, has the same eigenvalues of $B$.
Notice that, with $S$ diagonal positive definite, the products $BS$, $SB$ and $S^{\half}B S^{\half}$ are similar matrices,
\begin{gather*}
BS = S^{-1} (SB) S \sim SB
\\\text{and}\quad 
BS = S^{-\half} (S^{\half}B S^{\half}) S^{\half}\sim S^{\half}B S^{\half},  
\end{gather*}
which implies that they have the same eigenvalues. 
Moreover,
\begin{align*}
S^{\half}B S^{\half}
& = S^{\half} D^{\half} \sym{B} D^{-\half} S^{\half}
= D^{\half} (S^{\half} \sym{B} S^{\half}) D^{-\half}
\\
& \sim S^{\half} \sym{B} S^{\half},
\end{align*}
which implies that $S^{\half}B S^{\half}$ and $S^{\half} \sym{B} S^{\half}$ have the same eigenvalues.
Since $S^{\half} \sym{B} S^{\half}$ is a symmetric matrix, with $S^{\half}$ diagonal and nonsingular, it is possible to apply Theorem~\ref{theorem:Ostrowski}.
For all $k\in \{1,\mydots,n\}$ there is a positive real number $\theta_k \in [\min_i\{s_i\},\max_i\{s_i\}]$ such that 
\begin{gather*}
\lambda_{k}(S^{\half} \sym{B} S^{\half})= \theta_k \lambda_k(\sym{B}).
\end{gather*}
Therefore, from similarity, it follows that for all $k\in \{1,\mydots,n\}$ there exists $\theta_{k} \in [\min_i\{s_i\},\max_i\{s_i\}]$ such that 
\begin{gather*}
\lambda_{k}(S^{\half}B S^{\half})
=\lambda_{k}(SB) =\lambda_{k}(BS)
= \theta_k \lambda_k(B). \qedh
\end{gather*}%


\section{Proof of Theorem~\ref{thm:CT_SUB_summary}}
\label{proof:CT_SUB_Summary}
To improve readability, the proof of Theorem~\ref{thm:CT_SUB_summary} is divided as follows: in Section~\ref{sec:StabilityOrigin} we prove (i); in Sections~\ref{sec:ExistenceNonTrivialEq} and \ref{sec:StabilityNonTrivialEq}, \ref{sec:Uniqueness} we prove the existence (ii.1), stability (ii.2) and uniqueness (ii.3) part, respectively, of (ii).
The proof of (iii) is omitted since it is identical to the proof of (ii.1).

\subsection{Proof of Theorem~\ref{thm:CT_SUB_summary}(i)}
\label{sec:StabilityOrigin}
The condition for the existence of a unique equilibrium point for the system~\eqref{eqn:System_CT} can be rewritten in terms of the biggest eigenvalue of the normalized interaction matrix $H=\Delta^{-1}A= I -\Lnorm$.
Define the following symmetric matrix:
\begin{equation}
\sym{H} := \Delta^{-\half} A \Delta^{-\half} = \Delta^{\half} H \Delta^{-\half} \sim H,	\label{eqn:Hsym}
\end{equation}
By construction, $H$ and $\sym{H}$ have the same eigenvalues. 

\proof{}
First, notice that since $\lambda_i(\Lnorm)=1-\lambda_{n-i+1}(H)$ for all $i=1,\mydots,n$, then $\pi_1= \frac{1}{1-\lambda_1(\Lnorm)} = \frac{1}{\lambda_n(H)}$.

Let $V : \R^n \to \R_+$ be the Lyapunov function described by 
\begin{equation}
V(x)=\sum_{i=1}^{n}\,\int_{0}^{x_i} \psi_i(s)\,ds.
\label{eqn:Lyapunov}
\end{equation}
Since each function $\psi_i(\cdot)$ is monotonically increasing and $\psi_i(s)=0$ if and only if $s=0$, then $V(x)>0$ for all $x\in  \R^n\setminus \{0\}$ and $V(0)=0$. Moreover, $V(x)$ is radially unbounded.

From the assumptions \eqref{assumption:1psiOdd}, \eqref{assumption:2psiMonotone} and \eqref{assumption:4Sigmoidal}, we know that 
\begin{equation*}
x_i
\begin{cases} > \psi_i(x_i) ,&\text{if }x_i>0\; (\text{i.e., }\psi_i(x_i)>0)
	\\ < \psi_i(x_i) ,&\text{if }x_i<0\; (\text{i.e., }\psi_i(x_i)<0)
	\\ = 0 ,&\text{if }x_i=0,
\end{cases}
\end{equation*}
i.e., $\psi(x)^T \Delta x>\psi(x)^T \Delta \psi(x)>0$ since $x\ne 0$.
Hence, computing the derivative of $V$ along the trajectories gives
\begin{align*}
&\dot V(x)
= \psi(x)^T \dot x 
= \psi(x)^T [-\Delta x + \pi A \psi(x)]
\notag\\
&= - \psi(x)^T \Delta x + \psi(x)^T \Delta^\half (\pi \sym{H}) \Delta^\half \psi(x)
\notag\\
&< - \psi(x)^T \Delta^\half (I-\pi \sym{H}) \Delta^\half \psi(x)
\notag\\
&\le -\Bigl(1-\frac{\pi}{\pi_1}\Bigr) \psi(x)^T \Delta \psi(x)
\label{eqn:Vdot}
\end{align*}
Since $\pi\le \pi_1$, then $\dot V(x)<0$ for all $x\ne 0$, i.e., the origin is globally asymptotically stable, hence the unique equilibrium point for the system~\eqref{eqn:System_CT}.
\qed

\subsection{Proof of Theorem~\ref{thm:CT_SUB_summary}(ii.1): existence}
\label{sec:ExistenceNonTrivialEq}
Since $\lambda_1(\Lnorm)= 1-\lambda_n(H)$, the condition on $\pi$ can be rewritten as $\pi > \pi_1= \frac{1}{1-\lambda_1(\Lnorm)} = \frac{1}{\lambda_n(H)}$.
%
The proof follows \cite[Chapter I\S3]{GolubitskySchaeffer1985}.
The equilibrium points of the system~\eqref{eqn:System_CT} are solution of
\begin{equation}
\Phi(x,\pi) = -x + \pi H \psi(x) = 0.
\label{eqn:Proof_Eq}
\end{equation}
Let $J := \pde{\Phi}{x}(0,\pi_1) = -I + \pi_1 H$ be the Jacobian matrix at $(0,\pi_1)$, and let $v$ and $w$ be its left and right eigenvectors (such that $w^Tv = 1$), respectively, associated with the zero eigenvalue (notice that $\lambda_n(J) = -1+\pi_1 \lambda_n(H) =0$ by construction). 
By assumption $\lambda_n(H)$ is simple, which implies that $\dim(\ker{J}) = 1$ and that we can write $\ker{J} = \text{span}\{ v\}$ and $\range{J} =(\ker{J^T})^\perp = (\text{span}\{w\})^\perp$.
Let $E$ denote the projection of $\R^n$ onto $\range{J}$, $E = I- vw^T$, and $I-E = vw^T$ the projection onto $\ker{J}$.
The system of equations~\eqref{eqn:Proof_Eq} can be expanded as follows
\begin{subequations}
	\begin{align}
	E \,\Phi(x,\pi) &= 0	\label{eqn:Proof_Range}\\
	(I-E) \,\Phi(x,\pi) &= 0. \label{eqn:Proof_Ker}
	\end{align}
\end{subequations}
Split the vector $x$ accordingly, $x = yv+r$ with $y\in \R$, $yv \in \ker{J}$ and $r\in \range J$. 
From the implicit function theorem, solving \eqref{eqn:Proof_Range} for $r$ gives $r = R(yv,\pi)$, which we can substitute in \eqref{eqn:Proof_Ker} obtaining
\begin{align*}
0&= (I-E) \,\Phi(yv + R(yv,\pi),\pi)
\\&= (I-E) \,\bigl(-yv - R(yv,\pi) +\pi H  \psi(yv + R(yv,\pi))\bigr).
\end{align*}
Defining the center manifold $g:\R\times \R \to \R$ by 
\begin{equation*}
g(y) = w^T (I-E) \,\Phi(yv + R(yv,\pi),\pi)
\end{equation*}
we obtain that the zeros of $g$ are in one-to-one correspondence with the solutions of $\Phi(x,\pi)=0$.
We say that the system~\eqref{eqn:System_CT} undergoes a pitchfork bifurcation at $(0,\pi_1)$ if
\begin{equation*}
g = g_y = g_{yy} = g_\pi = 0,\quad g_{yyy} < 0,\quad  g_{\pi y} > 0,
\end{equation*}
where the subscript indicates partial derivative.

In what follows we will use the following notation: $\ppde{\psi}{x}(x) := \diag{\ppde{\psi_1}{x_1}(x_1),\dots,\ppde{\psi_n}{x_n}(x_n)}$ and $\ptde{\psi}{x}(x) = \diag{\ptde{\psi_1}{x_1}(x_1),\dots,\ptde{\psi_n}{x_n}(x_n)}$. Moreover, as $\Phi(x,\pi)$ is an odd function of $x$, we can neglect $R(\cdot)$ \cite[p.~33]{GolubitskySchaeffer1985}.
Calculations yield
\begin{subequations}
	\begin{align}
	\Phi(0,\pi_1) &= 0
	\label{eqn:Proof_Phi0}
	\\
	\Phi_y(0,\pi_1)
	&=\!\!\left.\left(-I+\pi H \pde{\psi}{x}(x) \right)\right|_{(0,\pi_1)} \!\! v 
	\notag\\&= (-I +\pi_1 H) v = 0
	\label{eqn:Proof_Phiy0}
	\\
	\Phi_{yy}(0,\pi_1)
	&=\left.\left(\pi H \;\pde{\psi}{x}\!\!\left(\pde{\psi}{x}(x) \,v\right)\right)\right|_{(0,\pi_1)} \!\! v 
	\notag
	\\&= \pi_1 H \;\ppde{\psi}{x}(0) \begin{bmatrix} v_1^2 \\ \vdots \\ v_n^2\end{bmatrix} = 0
	\label{eqn:Proof_Phiyy0}
	\\
	\Phi_{\pi}(0,\pi_1)
	&=H \psi(0) = 0 
	\label{eqn:Proof_Phipi0}
	\\
	\Phi_{yyy}(0,\pi_1)
	&=\left.\left( \pi H \;\pde{\psi}{x}\!\!\left(\pde{\psi}{x}\!\!\left(\pde{\psi}{x}(x) \,v\right)v\right)\right)\right|_{(0,\pi_1)} \!\! v 
	\notag
	\\&= \pi_1 H \;\ptde{\psi}{x}(0) \begin{bmatrix} v_1^3 \\ \vdots \\ v_n^3\end{bmatrix},
	\label{eqn:Proof_Phiyyy0}
	\\
	\Phi_{\pi y}(0,\pi_1)
	&=H \left.\pde{\psi}{x}(x)\right|_{(0,\pi_1)} \!\! v = H v = \lambda_n(H) v
	\label{eqn:Proof_Phipiy0}
	\end{align}
\end{subequations}
where $\psi(0)=0$, $\pde{\psi}{x}(0)=I$ and $\ppde{\psi}{x}(0)= 0$
directly follow from assumptions~\eqref{assumption:1psiOdd}, \eqref{assumption:2psiMonotone} and \eqref{assumption:4Sigmoidal}.

Equations~\eqref{eqn:Proof_Phi0}, \eqref{eqn:Proof_Phiy0}, \eqref{eqn:Proof_Phiyy0} and \eqref{eqn:Proof_Phipi0} yield respectively
\begin{align*}
g(0,\pi_1) = 0,\;
g_y(0,\pi_1) = 0,\;
g_{yy}(0,\pi_1) = 0,\;
g_\pi(0,\pi_1) = 0.
\end{align*}
The last two derivatives are given by
\begin{align*}
g_{\pi y}(0,\pi_1) &= \lambda_n(H) \,w^T v = \lambda_n(H)>0
\end{align*}
and
\begin{align*}
g_{yyy}(0,\pi_1)
& = \pi_1 w^T H \;\ptde{\psi}{x}(0) \begin{bmatrix} v_1^3 \\ \vdots \\ v_n^3\end{bmatrix}\\&
= \pi_1 \lambda_n(H) \sum_{i=1}^n w_i v_i^3 \;\ptde{\psi_i}{x_i}(0)<0.
\end{align*}
The last inequality holds since $\pi_1 \lambda_n(H) =1$, $w_i v_i\ge 0$ $\forall \,i$ and $\ptde{\psi_i}{x_i}(0)<0$ $\forall \,i$, which follows from assumption~\eqref{assumption:4Sigmoidal}.

We have shown that $(0,\pi_1)$ is a pitchfork bifurcation point for the system~\eqref{eqn:System_CT} and, as a consequence, that two new (nontrivial) equilibrium branches are created.
\qed

\subsection{Proof of Theorem~\ref{thm:CT_SUB_summary}(ii.2): stability}\label{sec:StabilityNonTrivialEq}
Before stating the proof we show that, when $\pi \in (\pi_1,\pi_2)$, all the nontrivial equilibrium points $x^*\ne 0 $ of the system~\eqref{eqn:System_CT} (if present) are locally asymptotically stable.
This lemma will be used in the proof of both the stability and uniqueness part of Theorem~\ref{thm:CT_SUB_summary}(ii).
Moreover, notice that the threshold values $\pi_1$ and $\pi_2$ can be rewritten as follows:
\begin{equation*}
\pi_1 =\frac{1}{\lambda_n(H)},\quad 
\pi_2 =\frac{1}{\lambda_{n-1}(H)}.
\end{equation*}%
\begin{lemma}\label{lemma:all_eq_are_stable}
Under the assumptions of Theorem~\ref{thm:CT_SUB_summary}(ii), when $\pi\in (\pi_1,\pi_2)$, if $x^\ast \ne 0$ is a (nontrivial) equilibrium point of the system~\eqref{eqn:System_CT} then it is locally asymptotically stable.
\end{lemma}
\proof{}
	Let $x^*\ne 0$ be an equilibrium point for the system~\eqref{eqn:System_CT}, 
	\begin{equation}
	x^* = \pi H \psi(x^*).
	\label{eqn:ASproof1}
	\end{equation}
	To prove that $x^*$ is locally asymptotically stable, consider the linearization around $x^*$,
	\begin{equation}
	\dot x = \Delta \Bigl(-I+ \pi H \pde{\psi}{x}(x^*)\Bigr)(x-x^*).
	\label{eqn:ASproof0}
	\end{equation}
	The equilibrium point $x^*$ is asymptotically stable for the system~\eqref{eqn:ASproof0}, and consequently locally asymptotically stable for the system~\eqref{eqn:System_CT}, if the matrix $\Delta \bigl(-I+ \pi H \pde{\psi}{x}(x^*)\bigr)$ is Hurwitz stable, i.e., its eigenvalues are strictly negative. Since $\Delta$ is diagonal and positive definite, this holds if and only if the matrix $-I+ \pi H \pde{\psi}{x}(x^*)$ is Hurwitz stable.
	The following proof shows that the largest eigenvalue of $-I+\pi H \pde{\psi}{x}(x^*)$ is strictly smaller than $0$, i.e., that $\pi \lambda_n\bigl(H \pde{\psi}{x}(x^*)\bigr)<1$. It is a two-steps proof, showing first that $\pi \lambda_n\bigl(H \pde{\psi}{x}(x^*)\bigr)\le 1$ and then by contradiction that $\pi \lambda_n\bigl(H \pde{\psi}{x}(x^*)\bigr)\ne 1$.
	
	\noindent\textbf{Step 1.}
	From the assumptions \eqref{assumption:1psiOdd}, \eqref{assumption:2psiMonotone} and \eqref{assumption:4Sigmoidal}, for all $i=1,\mydots,n$ it holds that 
	\begin{align*}
	\begin{cases}
	x_i^*>\psi(x_i^*)> \pde{\psi_i}{x_i}(x_i^*) \,x_i^*,\; &\text{if }x_i^*>0\\
	x_i^*<\psi(x_i^*)< \pde{\psi_i}{x_i}(x_i^*) \,x_i^*,\; &\text{if }x_i^*<0\\
	x_i^*=\psi(x_i^*)=0,\;\pde{\psi_i}{x_i}(x_i^*) =1,\; &\text{if }x_i^*=0.
	\end{cases}
	\end{align*}
	Therefore, $\exists \; \Xi(x^*)= \diag{\xi_1(x^*_1),\mydots,\xi_n(x^*_n)}\in \R^{n\times n}$ such that, $\forall \,i=1,\mydots,n$,
	\begin{gather*}
	\psi(x_i^*)= \xi_i(x_i^*) \cdot \pde{\psi_i}{x_i}(x_i^*) \cdot x_i^*,\\
	\begin{cases}
	\xi_i(x^*_i)> 1,\quad 0<\xi_i(x^*_i)\cdot \pde{\psi_i}{x_i}(x_i^*)<1,&\text{if }x_i^*\ne 0
	\\
	\xi_i(x^*_i)= 1,\quad \xi_i(x^*_i)\cdot \pde{\psi_i}{x_i}(x_i^*)=1,&\text{if }x_i^*= 0.
	\end{cases} 
	\end{gather*}
	This can be rewritten in compact form as
	\begin{equation}
	\psi(x^*)= \Xi(x^*)\cdot \pde{\psi}{x}(x^*) \cdot x^*,
	\label{eqn:ASproof3_1}	
	\end{equation}
	where
	\begin{equation}
	\diag{\Xi(x^*)}\ge \1,\;\;
	0< \diag{\Xi(x^*)\cdot \pde{\psi}{x}(x^*)}\le \1.
	\label{eqn:ASproof3_2}
	\end{equation}
	To simplify the notation, we will neglect the dependence from $x^*$ in what follows; moreover, we define $\dpsi := \pde{\psi}{x}(x^*)$.
	From \eqref{eqn:ASproof1} and \eqref{eqn:ASproof3_1}, it follows that
	\begin{equation}
	x^* = (\pi H \cdot\dpsi \cdot\Xi) \, x^*,
	\label{eqn:ASproof2}
	\end{equation}
	that is, $(1,x^*)$ is an eigenpair of $\pi H \,\dpsi \,\Xi$.
	Therefore, by Lemma \ref{lemma:DiagonalEquivalence} with $B=\pi H$ and $S =\dpsi \Xi$, and by \eqref{eqn:ASproof3_2}, it follows that for all $k\in \{1,\mydots,n\}$, $\exists \,\theta_k\in(0,1]$ such that $\lambda_{k}(\pi H \,\dpsi \,\Xi)= \theta_k \pi\lambda_k(H)$. In particular $ \exists \,\theta_{n}, \theta_{n-1} \in(0,1] $ such that
	\begin{gather}
	\lambda_{n}(\pi H \,\dpsi \,\Xi)= \theta_n \pi\lambda_n(H),
	\label{eqn:ASproof4_1}
	\\\lambda_{n-1}(\pi H \,\dpsi \,\Xi)= \theta_{n-1} \pi\lambda_{n-1}(H)
	< \theta_{n-1} \le 1,
	\label{eqn:ASproof4_2}	
	\end{gather}
	where $\pi\lambda_{n-1}(H)=\frac{\pi}{\pi_2}<1$ under the assumption $\pi\in (\pi_1,\pi_2)$.
	Then, \eqref{eqn:ASproof2} and \eqref{eqn:ASproof4_2} yield $1 = \lambda_{n}(\pi H \,\dpsi \,\Xi)$.
	
	Applying again Lemma~\ref{lemma:DiagonalEquivalence} with $B=\pi H\,\dpsi $ and $S=\Xi$, there exists $\theta\in [\min_i\{\xi_i\},\max_i\{\xi_i\}]$ such that
	\begin{equation*}
	1
	= \lambda_{n}(\pi H \,\dpsi \,\Xi)
	= \theta \pi \lambda_{n}(H \,\dpsi)
	\ge \pi \lambda_{n}(H \,\dpsi)>0,
	\end{equation*}
	since $\xi_i\ge 1$ for all $i$ implies $\theta \ge 1$. Hence, $\pi \lambda_{n}(H \,\dpsi)\le 1$.
	
	\noindent
	\textbf{Step 2.} Suppose by contradiction that $\pi \lambda_{n}(H \,\dpsi)=1$.
	First, we need to define the ``symmetric versions'' of the matrices $ H \,\dpsi $ and $ H \,\dpsi \, \Xi $.
	\begin{align*}
	\sym{[H\dpsi]}
	&:= (\dpsi)^{\half}\sym{H} (\dpsi)^{\half} 
	\\
	&\sim 
	(\Delta\dpsi)^{-\half} \!\cdot\! \Bigl((\dpsi)^{\half}\sym{H} (\dpsi)^{\half}\Bigr)\!\cdot\! (\Delta\dpsi)^{\half}
	\\&	=(\Delta^{-\half} \sym{H} \Delta^{\half})\,\dpsi
	= H \,\dpsi
	\end{align*}
	with $\sym{H}$ defined in \eqref{eqn:Hsym}, and
	\begin{align*}
	&\sym{[H\dpsi\Xi]}
	:=(\dpsi \,\Xi)^{\half}\sym{H} (\dpsi \,\Xi)^{\half}
	\\
	&\;\;\sim
	(\Delta\dpsi \,\Xi)^{-\half} \!\cdot\! \Bigl((\dpsi \,\Xi)^{\half}\sym{H} (\dpsi\, \Xi)^{\half}\Bigr)\!\cdot\! (\Delta\dpsi \,\Xi)^{\half}
	\\
	&\;\; = (\Delta^{-\half} \sym{H} \Delta^{\half})\, \dpsi\, \Xi
	= H \,\dpsi \, \Xi.
	\end{align*}
	Moreover, $\sym{[H\dpsi\Xi]} = \Xi^{\half}\,\sym{[H\dpsi]}\,\Xi^{\half}$ by construction.
	From similarity, $\pi \lambda_{n}(\sym{[H\dpsi]})=1$. Let $v$ be the right eigenvector of $\pi \,\sym{[H\dpsi]}$ associated with its largest eigenvalue $\pi \lambda_n(\sym{[H\dpsi]})=1$.
	From \eqref{eqn:ASproof2}, since $(1,x^*)$ is an eigenpair of $\pi H \,\dpsi \,\Xi$ (with $1$ being the largest eigenvalue, as proven previously), it follows that $(1,(\Delta\dpsi \,\Xi)^{\half} x^*)$ is an eigenpair of $\pi \,\sym{[H\dpsi\Xi]}$ (with $1$ being the largest eigenvalue).
	To summarize
	\begin{equation}
	\begin{cases}
	\sym{[H\dpsi\Xi]} = \Xi^{\half}\,\sym{[H\dpsi]}\,\Xi^{\half};
	\\
	\pi \lambda_n(\sym{[H\dpsi]})=1\;\;\text{and}\;\; \pi \sym{[H\dpsi]} \, v = v;
	\\
	\pi \lambda_n(\sym{[H\dpsi\Xi]}) = 1\;\;\text{and}\\
	\quad\quad \pi \sym{[H\dpsi\Xi]} \,(\Delta\dpsi \,\Xi)^{\half} x^* = (\Delta\dpsi \,\Xi)^{\half}x^*.
	\end{cases}
	\label{eqn:ASproof3}
	\end{equation}
	Applying Rayleigh's Theorem \cite[Thm 4.2.2]{HornJohnson2013} with $\pi \sym{[H\dpsi\Xi]}$, one obtains
	\begin{align}
	1
	&=\pi \lambda_n(\sym{[H\dpsi\Xi]}) 
	= \pi \;\max_{y\ne 0} \frac{y^T\;\sym{[H\dpsi\Xi]}\;y}{y^T y}
	\notag\\
	&= \pi \;\max_{y\ne 0} \frac{(y^T\; \Xi^{\half}) \cdot \sym{[H\dpsi]} \cdot (\Xi^{\half}\;y)}{y^T y}
	\notag\\
	&\quad \Downarrow \quad \text{with}\quad y = \Xi^{-\half} v\ne 0
	\notag\\
	&\ge \pi \;\frac{v^T\;\sym{[H\dpsi]} \;v}{v^T\Xi^{-1}v}
	= \frac{v^Tv}{v^T\Xi^{-1}v}.
	\label{eqn:ASproof5}
	\end{align}
	The inequality $v^Tv\le v^T\Xi^{-1}v$, which can be rewritten as
	\begin{equation*}
	\sum_{i=1}^n  v_i^2
	\le 
	\sum_{i=1}^n  \frac{v_i^2}{\xi_i}
	= \sum_{i:\, \xi_i=1} v_i^2
	+ \sum_{i:\, \xi_i>1} \frac{v_i^2}{\xi_i},	
	\end{equation*}
	holds, as equality, if and only if 
	\begin{equation}
	v_i=0 \quad \forall \,i \quad\text{s.t.}\quad \xi_i>1
	\;\; (\Leftrightarrow \; x_i^*\ne 0).
	\label{eqn:ASproof6}
	\end{equation}
	Hence, \eqref{eqn:ASproof5} can only hold as equality if \eqref{eqn:ASproof6} holds, which further implies that $\Xi^\half v = v$ and consequently that $(1,v)$ is an eigenpair of $\pi \sym{[H\dpsi\Xi]}$ (with $1$ being the largest eigenvalue). Indeed
	\begin{multline*}
	\pi \sym{[H\dpsi\Xi]} \,v
	= \pi \,\Xi^{\half}\,\sym{[H\dpsi]}\,\Xi^{\half} \,v
	\\
	= \pi \,\Xi^{\half}\,\sym{[H\dpsi]}\,v
	= \Xi^{\half}\,v = v.
	\end{multline*}

	Since $\lambda_n(\pi \sym{[H\dpsi\Xi]})=1$ is simple, as shown in \eqref{eqn:ASproof4_1} and \eqref{eqn:ASproof4_2}, it follows that $v$ should be equivalent to the corresponding right eigenvector of $\pi \sym{[H\dpsi\Xi]}$, i.e., $(\Delta\dpsi \,\Xi)^{\half} x^*$.
	This however yields a contradiction, since by \eqref{eqn:ASproof6} $v^T x^*=0$ (in particular $v_i=0$ for all $x_i^*\ne 0$ and $v_i\ne 0$ for at least one $i$ s.t. $x_i^*= 0$).
	
	To conclude, $\pi \lambda_{n}(H \,\dpsi)<1$ and the matrix $-I + \pi H \pde{\psi}{x}(x^*)$ is Hurwitz stable, which implies that the nontrivial equilibrium point $x^*$ is locally asymptotically stable.
	\qed

Finally, we are ready to prove Theorem~\ref{thm:CT_SUB_summary}(ii.2): stability.\\ 
\proof{Theorem~\ref{thm:CT_SUB_summary}(ii.2)}
The linearized system around the origin is $\dot x = \Delta(-I +\pi H) x$ (given assumption \eqref{assumption:2psiMonotone}), where $\Delta$ is positive definite and $-I+ \pi H$ has eigenvalues $\{ \pi \lambda_1(H)-1,\dots,\pi \lambda_n(H)-1\}$. When $\pi>\pi_1$, the matrix $-I+ \pi H$ has at least one positive eigenvalue, which proves the instability of the origin as equilibrium point of \eqref{eqn:System_CT}.

Instead, let $x^\ast\ne 0 $ be an equilibrium point of the system~\ref{eqn:System_CT}, whose existence is shown in Theorem~\ref{thm:CT_SUB_summary}(ii.1).
To prove its stability we can apply Lemma~\ref{lemma:all_eq_are_stable}, which shows that if $x^*\ne 0$ is an equilibrium point of \eqref{eqn:System_CT} then it must be locally asymptotically stable.
\qed

\subsection{Proof of Theorem~\ref{thm:CT_SUB_summary}(ii.3): uniqueness}
\label{sec:Uniqueness}
Let $f(x,\pi)=\Delta[-x+\pi H \psi(x)]$.
We can divide the proof into two parts. 
(i) We prove that in a neighborhood of $\pi_1$ the system admits only 3 equilibria. 
(ii) We prove, by contradiction, that there are no bifurcation points for $\pi \in (\pi_1,\pi_2)$. We use Lemma~\ref{lemma:all_eq_are_stable} which proves that \textit{all} the nontrivial equilibrium points $\bar x$ of the system~\eqref{eqn:System_CT} are locally asymptotically stable, hence isolated and with each matrix $\pde{f}{x}(\bar{x},\bar{\pi})$ Hurwitz.

\textbf{(i).}
The existence is shown in Theorem~\ref{thm:CT_SUB_summary}(ii.1), where it is also proven that the bifurcation is a pitchfork. This means that in a neighborhood of $\pi_1$ the system~\eqref{eqn:System_CT} admits exactly three equilibrium points: the origin and two nontrivial equilibrium points, $\pm x^*\ne 0$.

\textbf{(ii).}
The necessary condition for an equilibrium point $(\bar{x}, \bar{\pi})$ (where $\bar \pi \in (\pi_1,\pi_2)$)  to be a bifurcation point is that the Jacobian $\pde{f}{x}(\bar{x},\bar{\pi}) = \Delta [- I +\bar{\pi} H \pde{\psi}{x}(\bar x) ]$ is not invertible (i.e., there is an $i\in \{1,\mydots,n\}$ such that $\bar{\pi} \lambda_i(H \pde{\psi}{x}(\bar x))=1$).
Suppose by contradiction that, for $\bar \pi \in (\pi_1,\pi_2)$, $\bar{x}$ is an equilibrium point of the system~\eqref{eqn:System_CT}, i.e., $\bar x = \bar \pi H \psi(\bar x)$, and a bifurcation point, i.e., $\exists \; i$ s.t. $\bar{\pi} \lambda_i(H \pde{\psi}{x}(\bar x))=1$.
However, Lemma~\ref{lemma:all_eq_are_stable} shows that if $\bar x\ne 0$ is an equilibrium point of \eqref{eqn:System_CT} then it must be locally asymptotically stable., i.e., $\bar \pi \lambda_n(H \pde{\psi}{x}(\bar x))<1$. Moreover, if $\bar x = 0$, $\bar \pi \lambda_{n-1}(H)<1<\bar \pi \lambda_{n}(H)$ for $\bar \pi \in (\pi_1,\pi_2)$.
Hence, $\bar x$ cannot be a bifurcation point.

To conclude, we know that the system~\eqref{eqn:System_CT} admits three equilibria ($0$, $x^*$, $-x^*$) and that it cannot bifurcate further from them for values of $\pi \in (\pi_1,\pi_2)$. 
Hence, the only possible equilibrium points for the system are the origin and those originated from the first bifurcation at $\pi =\pi_1$.
\qed

\section{Proof of Lemma~\ref{lemma:lambda2}}
\label{proof:lambda2}
It is useful to first reformulate structural balance of $\G$ (see Section~\ref{sec:Preliminaries}) in terms of the matrix $\sym{H}$ defined in \eqref{eqn:Hsym}, which is congruent to $A$ and similar to $H$: $\G$ (connected) is structurally balanced if and only if there exists a signature matrix $S=\diag{s_1,\mydots,s_n}$ with diagonal entries $s_i=\pm 1$ ($i=1,\mydots,n$) such that $S\sym{H}S$ is nonnegative.

We are now ready to prove Lemma~\ref{lemma:lambda2}. In what follows we say that a symmetric matrix is ``special'' if it has zero diagonal entries, and ``elliptic'' if it has exactly one and simple positive eigenvalue. This notation is from \cite{Fiedler1994}, whose results are used in the proof.
Notice that both $A$ and $\sym{H}$ are special matrices.\\
\begin{proof}{}
	Let $\G$ be structurally unbalanced.
	
	Assume by contradiction that $\lambda_2(\Lnorm)> 1$, which is equivalent to $\lambda_{n-1}(H)=\lambda_{n-1}(\sym{H})<0$.
	This means that $\sym{H}$ is special, nonsingular and elliptic.
	Then, all the off-diagonal entries of $\sym{H}$ are different from zero (see \cite[Corollary~2.7]{Fiedler1994}) and hence there exists a diagonal (signature) matrix $S = \diag{s_1,\mydots, s_n}$ with $s_i=\pm 1$ s.t. $S\sym{H}S$ is nonnegative (see \cite[Thm~2.5]{Fiedler1994}). Therefore, $\G$ is structurally balanced and we obtain a contradiction.
		
	Assume by contradiction that $\lambda_2(\Lnorm)=1$, which is equivalent to $\lambda_{n-1}(H)=\lambda_{n-1}(\sym{H})=0$. 
	This means that $\sym{H}$ is special, singular and elliptic.
	Let $r= \text{rank}(\sym{H})<n$ and observe that $\sym{H}$ cannot have zero rows or columns since $\G$ is connected.
	From \cite[Thm~2.9]{Fiedler1994}, there exist a permutation matrix $P\in \R^{n\times n}$ and an integer $t\in \{r,\mydots,n\}$ s.t. 
	\begin{equation*}
		P\sym{H} P^T = D {\sym{H}}_0 D^T,
	\end{equation*}
	where ${\sym{H}}_0\in \R^{t\times t}$ is elliptic, special and nonnegative with rank $r$ and $D= d_1\oplus d_2 \oplus \mydots \oplus d_t \in \R^{n\times t}$ (here $\oplus$ indicates the direct sum) is a block matrix 
	where 
	each $d_i \in \R^{n_i}$ (with $\sum_{i=1}^t n_i =n$) is a unit vector (i.e., $\norm{d_i}_2=1$) with all elements different from zero (i.e., $\abs{d_i}>0$).
	Define the signature (block) matrix $S:= S_1\oplus \mydots \oplus S_t\in \R^{n\times n}$, where each $S_i:=\diag{\sign{d_i}}\in \R^{n_i\times n_i}$, $i=1,\mydots,t$, is a signature matrix.
	Then, $P\sym{H} P^T$ can be rewritten as
	\begin{align*}
		P\sym{H} P^T
		= S \abs{D} {\sym{H}}_0 \abs{D^T} S.
	\end{align*} 
	It follows that $S (P\sym{H} P^T) S$ is nonnegative, and hence, that $\bar{S} \sym{H} \bar{S} $ is nonnegative, with $\bar S = P^T S P $ still a signature matrix. Therefore, $\G$ is structurally balanced and we obtain a contradiction.
\end{proof}


\section{Proof of Theorem~\ref{thm:CT_BoundedEq}}\label{sec:BoundedTraj}
In the following proofs we use the notation $S_x:= \diag{\sign{x}}$ where $x\in \R^n$ and the signum function is defined as $\sign{y}=1$ if $y\le 0$ or $\sign{y}=-1$ if $y< 0$, where $y\in \R$.
\begin{remark}\label{remark:frustr_H}
The frustration of the network $\G$ is defined in equation \eqref{eqn:frustration}, which can be rewritten as follows:
\begin{align*}
2\frustration 
&= \min_{\substack{S=\diag{s_1,\mydots,s_n}\\s_i=\pm 1\; \forall i}} \1^T (\abs{H}-S H S)\1 
\\
&= n - \max_{\substack{S=\diag{s_1,\mydots,s_n}\\s_i=\pm 1\; \forall i}} \1^T S H S\1
\\
&= n - \max_{\substack{S=\diag{s_1,\mydots,s_n}\\s_i=\pm 1\; \forall i}} \1^T S \left(\frac{H+H^T}{2}\right) S\1.
\end{align*}
\end{remark}
\begin{remark}
\label{remark:ustar}
Let 
\begin{equation}
u^\ast=S_{u^\ast}\1\;\,\text{where}\,\; S_{u^\ast} = \argmax_{\substack{S=\diag{s_1,\mydots,s_n}\\s_i=\pm 1\; \forall i}} \1^T S H S\1,
\label{eqn:ustarSymm}
\end{equation}
that is, $2\frustration = n- \1^T S_{u^\ast} H S_{u^\ast} \1$.
From the results on (symmetric) Hopfield neural networks (see \cite{Hopfield1984,HopfieldTank1985}) we know that the vector $u^\ast$ satisfies $u^\ast = \sign{\bigl(\frac{H+H^T}{2}\bigr) u^\ast }$. If $H$ is symmetric, then $u^\ast = \sign{H^T u^\ast }$, meaning that the vector $S_{u^\ast} \sign{H^T S_{u^\ast} \1 }= S_{u^\ast} \sign{H S_{u^\ast} \1}$ has all strictly positive components (equal to $1$).	
\end{remark}

\subsection{Proof of Theorem~\ref{thm:CT_BoundedEq}(i)}
The proof can be divided into three steps. First, we show that if $x^\ast$ is an equilibrium point of the system~\eqref{eqn:System_CT} and $H$ is symmetric, then $\norm{x^\ast}_1 \le \pi (n-2\frustration)$.
Then, we show that if $H$ is not symmetric and $x^\ast$ is such that $x_i^\ast \ne 0$ for all $i$ (or, $\abs{x^\ast}>0$) we can apply the same reasoning to prove that $\norm{x^\ast}_1 \le \pi (n-2\frustration)$.
Finally, we complete the proof and show that each equilibrium point $x^\ast$ of the system~\eqref{eqn:System_CT} (without assuming that $H$ is symmetric) satisfies the inequality $\norm{x^\ast}_1 \le \pi (n-2\frustration)$.

\textbf{Step 1.} We first consider the particular case of the matrix $H = \Delta^{-1} A $ being symmetric, that is, $\Delta = \delta I$.

Let $x^\ast$ be an equilibrium point of the system~\eqref{eqn:System_CT}, that is, $x^\ast = \pi H \psi(x^\ast)$ and let $S_{x^\ast}$ be its signature, i.e., $\abs{x^\ast} = S_{x^\ast} x^\ast$. 
It follows that 
\begin{equation*}
\abs{x^\ast} = \pi\,S_{x^\ast} H \psi(x^\ast) = \pi\,\abs{H \psi(x^\ast)} 
= \pi\, S_{x^\ast} H S_{x^\ast} \abs{\psi(x^\ast)}.
\end{equation*} 
Observe that $S_{H \psi(x^\ast)}= S_{x^\ast} = S_{\psi(x^\ast)}$.
Then 
\begin{equation*}
\frac{\norm{x^\ast}_1}{\pi} = \frac{\1^T\abs{x^\ast}}{\pi} = \1^T \abs{H \psi(x^\ast)} = \1^T S_{x^\ast} H S_{x^\ast} \abs{\psi(x^\ast)}
\end{equation*}
and
\begin{align*}
\max_{\substack{x^\ast \in \R^n \text{s.t.}\\x^\ast = \pi H \psi(x^\ast)}} \frac{\norm{x^\ast}_1}{\pi}
&=\max_{\substack{x^\ast \in \R^n \text{s.t.}\\x^\ast = \pi H \psi(x^\ast)}} 
\1^T \abs{H \psi(x^\ast)}
\\&
\le \max_{\substack{u \in \R^n\text{s.t.}\\S_{H u} = S_{u},\;-\1 \le u \le \1}} \1^T \abs{H u}.
\end{align*}
because the constraint $x^\ast = \pi H \psi(x^\ast)$ (i.e., $x^\ast$ is an equilibrium point) implies the constraint $S_{H \psi(x^\ast)} =  S_{\psi(x^\ast)}$, and $ \abs{\psi_i(x_i)} \le 1 $ $\forall$ $x_i \in \R$.
Then
\begin{subequations}
\begin{align}
\max_{\substack{x^\ast \in \R^n \text{s.t. }\\x^\ast = \pi H \psi(x^\ast)}} &\frac{\norm{x^\ast}_1}{\pi}
\le \max_{\substack{u \in \R^n\text{s.t. }\\S_{H u} = S_{u},\;-\1 \le u \le \1}} \1^T \abs{H u}
\notag
\\
&= \max_{u \in \R^n\;\text{s.t. }-\1 \le u \le \1} \1^T S_u H S_u \abs{u}
\label{eqn:MaxProbSymm_1}
\\
&= \max_{\substack{S_u=\diag{s_{u,1},\mydots,s_{u,n}}\\s_{u,i} = \pm 1\;\forall i}} \1^T S_u H S_u \1
\label{eqn:MaxProbSymm_2}
\\
& = n-2\frustration.\notag
\end{align}
\end{subequations}
Notice that $\1^T S_u H S_u \abs{u} \le \abs{\1^T S_{u} H S_{u}}\1$ for all $u\in \R^n$ s.t. $\abs{u}\le \1$. The equality between \eqref{eqn:MaxProbSymm_1} and \eqref{eqn:MaxProbSymm_2} means that the maxima are obtained when $u$ lies in the corners of the hypercube $\abs{u}\le\1$ (i.e., $\abs{u}=\1$).
In particular, $u^\ast$ defined in \eqref{eqn:ustarSymm} is a solution of this maximization problem since it is feasible and $\sign{\1^T S_{u^\ast} H S_{u^\ast} }= \sign{\1^T S_{u^\ast} H}S_{u^\ast} = \1^T \ge 0$, meaning that $\1^T S_{u^\ast} H S_{u^\ast} = \abs{\1^T S_{u^\ast} H S_{u^\ast} }$. 

\textbf{Step 2.} Let $x^\ast$ be an equilibrium point of the system~\eqref{eqn:System_CT} and assume that $x_i^\ast \ne 0$ for all $i$. In this step we do not assume the symmetry of the matrix $H$.
Following the reasoning of \textbf{Step 1}, and by adding the additional constraint $\abs{u}>0$ (which comes from $x^\ast_i\ne 0$ for all $i$), it is still possible to prove (see below) that the maxima are obtained when $u$ lies in the corners of the hypercube $\abs{u}\le \1$ (which yields the equivalence between \eqref{eqn:MaxProbSymm_1} and \eqref{eqn:MaxProbSymm_2}).
This is equivalent to show that the maxima are obtained when $u$ is s.t. $\1^T S_u H S_u \ge 0$. Let 
\begin{equation*}
	\tilde{u} = \argmax_{u \in \R^n\;\text{s.t. }-\1 \le u \le \1,\,\abs{u} >0} \1^T S_u H S_u \abs{u}
\end{equation*}
and $v^T=\1^T S_{\tilde{u}} H S_{\tilde{u}}$.
Suppose, by contradiction, that $\exists\,j$ s.t. $v_j<0$ (and $\tilde u_j \ne 0$, since $\abs{\tilde u}>0$).
Define $\bar u$ s.t. $S_{\bar u} = S_{\tilde u}$ and
\begin{equation*}
	\abs{\bar u_j} = \begin{cases}
		1, &\text{if }v_j\ge 0
		\\ \varepsilon \abs{\tilde u_j}, & \text{if }v_j<0
	\end{cases}
\end{equation*} 
with $0<\varepsilon<1$ (which means that $\abs{\bar u}\le \1$ and $\abs{\bar u}> 0$, i.e., $\bar u$ is a feasible point of the maximization problem). Then
\begin{align*}
	\1^T S_{\bar u}& H S_{\bar u} \abs{\bar u} 
	= v^T \abs{\bar u} 
	= \sum_{j: v_j\ge 0} v_j \abs{\bar u_j} + \sum_{i: v_j< 0} v_j \abs{\bar u_j}
	\\
	&= \sum_{j: v_j\ge 0} v_j + \sum_{j: v_j< 0} v_j \varepsilon \abs{\tilde u_j}
	\\
	&\ge \sum_{j: v_j\ge 0} v_j \abs{\tilde u_j}+ \sum_{j: v_j< 0} v_j \varepsilon \abs{\tilde u_j}
	\\
	&> \sum_{j: v_j\ge 0} v_j \abs{\tilde u_j} + \sum_{i: v_j< 0} v_j \abs{\tilde u_j}
	= \1^T S_{\tilde u} H S_{\tilde u} \abs{\tilde u},
\end{align*}
which implies a contradiction.

\textbf{Step 3.} Finally, we want to extend the idea presented in the previous steps to show, without imposing any constraint on $H$ or $x^\ast$ (except for being an equilibrium point of the system~\eqref{eqn:System_CT}), that $\norm{x^\ast}_1\le \pi (n-2\frustration)$.
In this case, assume that $x_i^\ast=0$ for some $i$: let $n_0$ (resp., $n-n_0$) be the number of zero (resp., nonzero) components of $x_i^\ast$ and let $P$ be a permutation matrix s.t. $P x^\ast = \begin{bmatrix} x^\ast_{nz} \\ 0\end{bmatrix}$, where $x^\ast_{nz}\in \R^{n-n_0}$ and $\abs{x^\ast_{nz}}>0$. Let $P \psi(x^\ast) = \begin{bmatrix} \psi(x^\ast_{nz}) \\ 0\end{bmatrix}$ and $P H P^T = \begin{bmatrix} H_{nz} & \star\\ \star & \star\end{bmatrix}$. Then, $\norm{x^\ast}_1 = \norm{x^\ast_{nz}}_1$, $x^\ast_{nz} = \pi H_{nz}\psi(x^\ast_{nz})$ and (following the reasoning of \textbf{Step 2})
\begin{equation*}
\frac{\norm{x^\ast}_1}{\pi} 
\le \max_{\substack{S_{nz}=\diag{s_1,\mydots,s_{n-n_0}}\\s_i=\pm 1\; \forall i}} \1_{n-n_0}^T S_{nz}  H_{nz} S_{nz} \1_{n-n_0}.
\end{equation*}
Notice that for all signature matrices $S_{nz}=\diag{s_1,\mydots,s_{n-n_0}} \in \R^{n-n_0,n-n_0}$ with $s_i=\pm 1$ $\forall i$, the following inequality holds:
\begin{align*}
	\1_{n-n_0}^T &S_{nz}  H_{nz} S_{nz} \1_{n-n_0} 
	\\&= 
	\Bigl(\begin{bmatrix} \1_{n-n_0}^T & 0_{n_0}^T \end{bmatrix} P\Bigr)\cdotp \Bigl(P^T\begin{bmatrix} S_{nz} & 0 \\0 & I_{n_0} \end{bmatrix} P\Bigr) \\&\qquad \cdotp H \cdotp \Bigl(P^T \begin{bmatrix} S_{nz} & 0 \\0 & I_{n_0}  \end{bmatrix}P\Bigr)\cdotp \Bigl(P^T \begin{bmatrix} \1_{n-n_0} \\ 0_{n_0} \end{bmatrix}\Bigr)
	\\
	& \le \max_{u \in \R^n \text{s.t. }\abs{u}\le 1} u^T H u
	= n-2\frustration.
\end{align*}
Summarizing, we have shown that $\norm{x^\ast}_1 \le \pi (n-2\frustration) $.
\qed

\subsection{Proof of Theorem~\ref{thm:CT_BoundedEq}(ii)}
From Theorem~\ref{thm:CT_BoundedEq}(i) we know that the set
\begin{equation*}
\Omega_{\frustration} = \{x\in \R^n: \norm{x}_1\le \pi (n-2\frustration)\}
\end{equation*}
contains all equilibria of the system~\eqref{eqn:System_CT}. We now want to prove that $\Omega_{\frustration}$ is attractive when the signed normalized Laplacian $\Lnorm$ or, equivalently, the normalized interaction matrix $H = \Delta^{-1} A =I-\Lnorm$, is symmetric (i.e., $\Delta = \delta I$).

Let $V : \R^n \to \R_+$ be the Lyapunov function described by 
\begin{equation}
V(x)= \begin{cases}
\frac{1}{\delta}\bigl(\norm{x}_1 - \pi (n-2\frustration)\bigr), & x \notin \Omega_{\frustration}
\\
0, & x \in \Omega_{\frustration} \end{cases}.
\label{eqn:LyapunovSymm_BT}
\end{equation}
Since $\Delta$ is positive definite then $V(x)>0$ for all $x \notin \Omega_{\frustration}$. Moreover, $V(x)$ is radially unbounded.
Let $S_x:= \diag{\sign{x}}$ (observe that $S_{x}= S_{\psi(x)}$).
The upper Dini derivative of $V$ along the trajectories \eqref{eqn:System_CT} (with $\Delta = \delta I$) gives
\begin{align*}
d^+ V(x)
&= \limsup_{s\to 0^+} \frac{V(x+s \dot x) - V(x)}{s}
\\
&= \frac{1}{\delta} \limsup_{s\to 0^+} \frac{\sum_i \abs{x_i+s \dot x_i}-\abs{x_i} }{s}
\\
&= \frac{1}{\delta} \sum_i d^+\abs{x_i}
= \frac{1}{\delta} \sum_{i=1}^{n} \sign{x_i}\dot x_i
\\
&= \1^T \Delta^{-1} S_x \dot x
= \1^T S_x [-x+\pi H \psi(x)]
\\
&= -\norm{x}_1 + \pi \1^T S_x H \psi(x)
\\
&= -\norm{x}_1 + \pi \1^T S_x H S_x \abs{\psi(x)},
\\
&\le -\norm{x}_1 + \pi  \max_{u \in \R^n \text{ s.t. }\abs{u}\le \1} 1^T S_u H u
\end{align*}
Again, the intuition is that $\1^T S_u H u = \1^T S_u S_{H u} \abs{H u} \le \1^T \abs{H u}$, which means that the maxima of $\1^T S_u H u$ are obtained when $S_{H u} = S_u$.
Hence, following the reasoning of the proof of Theorem~\ref{thm:CT_BoundedEq}, we conclude that for all $x\notin \Omega_{\frustration}$
\begin{equation*}
d^+ V(x)\le -\norm{x}_1 + \pi (n-2\frustration)< 0.\qedh
\end{equation*}


\section{Proof of Proposition~\ref{prop:fpi1pi2}}
\label{sec:prop:fpi1pi2}
To show that the bound~\eqref{eqn:fpi1pi2} holds, observe first that $1\le \pi_1\le \pi_2$ holds trivially since $\frac{1}{\pi_1}=1-\lambda_1(\Lnorm)\le 1$ and $\pi_1 = \frac{1}{1-\lambda_1(\Lnorm)}\le \frac{1}{1-\lambda_{2}(\Lnorm)}=\pi_2$. 
Therefore, it remains to show that $\pi_1 \le \frac{n}{n-2\frustration}$.
It is easier to write the rest of the proof in terms of the normalized interaction matrix $H=I-\Lnorm=\Delta^{-1}A$; then, $\pi_1 = \frac{1}{1-\lambda_1(\Lnorm)}= \frac{1}{\lambda_n(H)}$.
Let $\Sbest$ be the matrix yielding the minimum value of energy in \eqref{eqn:frustration}, that is, $2\frustration = n - \1^T \Sbest H \Sbest \1$ (Remark~\ref{remark:frustr_H} in the proof of Theorem~\ref{thm:CT_BoundedEq} shows how the frustration of the network $\G$ can be rewritten in terms of $H$).
Since $\Lnorm$ is symmetric then $H$ is also symmetric and using the Rayleigh's Theorem \cite[Thm 4.2.2]{HornJohnson2013}, it follows that
\begin{equation*}
	\frac{1}{\pi_1} =\lambda_n(H) 
	= \max_{v\in \R^n} \frac{v^T H v}{v^T v}
	\ge \frac{\1^T \Sbest H \Sbest \1}{\1^T \1}
	= \frac{n-2\frustration}{n},
\end{equation*}
which implies \eqref{eqn:fpi1pi2}.
\qed


\section{Proof of Proposition~\ref{prop:fpi1pi1d}}\label{sec:DT_fpi1pi1d}
The first part of the proof holds for both structurally balanced and structurally unbalanced (connected) graphs $\G$ with normalized signed Laplacian $\Lnorm$. 

The condition $\pid>\pi_1$ holds if at $\pi_1$ the biggest eigenvalue of $L_{\pi_1}$ is smaller than $\frac{2}{\varepsilon}$, that is,
$\lambda_n(L_{\pi_1})<\frac{2}{\varepsilon}= \lambda_n(L_{\pid})$.
Since $L_{\pi_1}= \Delta(I-\pi_1 (I-\Lnorm))$, from Lemma~\ref{lemma:DiagonalEquivalence} it follows that $\exists\, \theta\in \bigl[\min_i \delta_i,\max_i \delta_i\bigr]$ such that 
\begin{equation*}
	\lambda_n(L_{\pi_1}) 
	= \theta \lambda_n(I-\pi_1 (I-\Lnorm))
	= \theta (1-\pi_1 (1-\lambda_n(\Lnorm))).
\end{equation*}
Let $\varepsilon \max_i \delta_i <1$. Observe that the smallest and largest eigenvalues of $\Lnorm$ always satisfy $\lambda_1(\Lnorm)<1<\lambda_n(\Lnorm)\le 2$ (see, e.g., \cite{LiLi2009}). Then $1-\lambda_n(\Lnorm)<0$ and
\begin{align*}
	\lambda_n&(L_{\pi_1})
	= \theta \left(1+\pi_1 (\lambda_n(\Lnorm)-1)\right)
	\\&\le \max_i \delta_i (1+\pi_1 (\lambda_n(\Lnorm)-1))
	< \frac{1}{\varepsilon} (1+\pi_1 (\lambda_n(\Lnorm)-1)).
\end{align*}
Therefore, showing that $\pi_1(\lambda_n(\Lnorm)-1)\le 1$ is sufficient to conclude that $\lambda_n(L_{\pi_1})<\frac{2}{\varepsilon}$ and hence that $\pid>\pi_1$. We now need to treat the cases (i) and (ii) separately.

Assume that (i) holds, i.e., $\lambda_1(\Lnorm) =0$. Then, since $\lambda_n(\Lnorm)\le 2$, $\pi_1(\lambda_n(\Lnorm)-1)=\lambda_n(\Lnorm)-1\le 1$.

Assume that (ii) holds, i.e., $\lambda_1(\Lnorm)<2-\lambda_n(\Lnorm)$. Then, $\pi_1 (\lambda_n(\Lnorm)-1)<\pi_1 (1-\lambda_1(\Lnorm))=1$.

Hence, the condition $\pid>\pi_1$ holds. \qed

\section{Proof of Lemma~\ref{lemma:DT_Origin} and Lemma~\ref{lemma:DT_Origin_Oscill}}
\label{sec:DT_NecCondNontrivialEq}
\proof{Lemma~\ref{lemma:DT_Origin}}
	Assume, by contradiction, that $\pi\le\pi_1$ which implies that $\Lpi$ is positive semidefinite (positive definite if $\pi<\pi_1$), see Remark~\ref{remark:eig_Lpi}.
	Let $x^\ast\ne 0$ be an equilibrium point of \eqref{eqn:System_DT}, that is, $\Delta x^\ast = \pi A \psi(x^\ast) = \Delta \psi(x^\ast) - \Lpi \psi(x^\ast)$.
	Then,
	\begin{equation*}
	0< \psi(x^\ast)^T \Delta (x^\ast -\psi(x^\ast)) = - \psi(x^\ast)^T \Lpi \psi(x^\ast)
	\le 0,
	\end{equation*}
	which leads to a contradiction. Hence, $\pi>\pi_1$.\qed

\proof{Lemma~\ref{lemma:DT_Origin_Oscill}}
	Assume, by contradiction, that $\pi\le\pid$ which implies that $\varepsilon\Lpi-2 I$ is negative semidefinite (negative definite if $\pi<\pid$), see Remark~\ref{remark:eig_Lpi}.
	Assume that the system~\eqref{eqn:System_DT} admits a period-2 limit cycle: $\exists\; K>0$ such that $x_{k+2} = \xk \ne 0$ for all $k\ge K$, that is
	\begin{equation*}
	\begin{cases}
	\xk = (I-\varepsilon \Delta) \xkp + \varepsilon\pi A \psi(\xkp)
	\\
	\xkp = (I-\varepsilon \Delta) \xk + \varepsilon\pi A \psi(\xk),
	\end{cases}
	\end{equation*}
	which implies that
	\begin{equation*}
	0 = (2I-\varepsilon \Delta) (\xkp-\xk) + \varepsilon\pi A (\psi(\xkp)-\psi(\xk)).
	\end{equation*}	
	Then,
	\begin{align*}
	&0 
	= (\psi(\xkp)-\psi(\xk))^T (2I-\varepsilon \Delta) (\xkp-\xk) 
	\\&\quad\;+ \varepsilon\pi (\psi(\xkp)-\psi(\xk))^T A (\psi(\xkp)-\psi(\xk))
	\\
	&>(\psi(\xkp)-\psi(\xk))^T (2I-\varepsilon \Delta +\varepsilon\pi A) (\psi(\xkp)-\psi(\xk))
	\\
	&= (\psi(\xkp)-\psi(\xk))^T (2I-\varepsilon \Lpi) (\psi(\xkp)-\psi(\xk))
	\ge 0
	\end{align*}
	which leads to a contradiction since $\psi(\xkp)-\psi(\xk)\ne 0$. Hence, $\pi>\pid$.
	Observe that the first inequality holds since $\varepsilon \max_i \delta_i <2$ and each nonlinearity is monotonically increasing (i.e., if $x_i(k+1)-x_i(k)\ge 0$ then $\psi_i(x_i(k+1))-\psi_i(x_i(k))\ge 0$ for all $i$) and Lipschitz with constant $1$ (i.e., $|x_i(k+1)-x_i(k)|\ge |\psi_i(x_i(k+1))-\psi_i(x_i(k))|$ for all $i$).
	\qed


\section{Proof of Theorem~\ref{thm:DT_summary}}\label{proof:DT_summary}
To improve readability, the proof of Theorem~\ref{thm:DT_summary} is divided as follows: in Section~\ref{proof:DT_StabilityOrigin} we prove (i) and in Section~\ref{proof:DT_PeriodDoubling_Pitchfork} we prove (ii).
	
\subsection{Proof of Theorem~\ref{thm:DT_summary}(i)}
\label{proof:DT_StabilityOrigin}
It is useful to first introduce the following two lemmas.
\begin{lemma}\label{lemma:DT_invconvex}
	Consider the function $f: [-1,1] \to \R$ defined as
	\begin{equation*}
	f(y) = \int_{0}^{y} \psi^{-1}(s) \der s - \frac{1}{2} y^2.
	\end{equation*}
	where $\psi:\R \to [-1,1]$ is a nonlinear function satisfying the properties \eqref{assumption:1psiOdd}$\div$\eqref{assumption:4Sigmoidal} and $\psi^{-1}$ indicates the inverse function.
	Then, $f(y)> 0$ for all $y\ne 0$ and $f(0)=0$.
\end{lemma}
\proof{}
	Form the definition of $f$ it follows that $f(0)=0$ and $f'(y)= \psi^{-1}(y)-y$.
	Since $f'(0)=0$, a sufficient condition for $f(y)> 0$ to hold for all $y\ne 0$ is that $f$ is a convex function.
	
	To prove that $f$ is convex we compute the second derivative:
	\begin{equation*}
		f''(y)= \frac{1}{\psi'(\psi^{-1}(y))} -1 
		\begin{cases}
		> 1-1 = 0, & y\ne 0\\
		= 0, & y=0
		\end{cases}
	\end{equation*} 
	since $0<\psi'(y)<1$ $\forall \,y\ne 0$ and $\psi'(0)=1$.
	It follows that $f$ is convex, which implies that $f(y)\ge 0$ $\forall \,y\in [-1,1]$. \qed

\begin{lemma}[Taylor expansion]\label{lemma:DT_invtaylor}
Consider the function $g : [-1,1] \to \R_+$ given by $g(y) = \int_{0}^{y} \psi^{-1}(s) \der s$, where $\psi:\R \to [-1,1]$ is a nonlinear function satisfying the properties \eqref{assumption:1psiOdd}$\div$\eqref{assumption:4Sigmoidal}.
Expanding $g$ around $y_0$ yields
\begin{align*}
\int_{0}^{y} \psi^{-1}(s) \der s
= &\int_{0}^{y_0} \psi^{-1}(s) \der s
\\&+ (y-y_0) \psi^{-1}(y_0) + \frac{(y-y_0)^2}{2} g''(z)
\end{align*}
where $z \in [y_0,y]$ and $g''(z):=\frac{1}{\psi'(\psi^{-1}(z))}\ge 1$.
\end{lemma}

Now we are ready to prove Theorem~\ref{thm:DT_summary}(i); the proof follows the work \cite{Zhaoetal2002}.\\
\proof{}
	Let $\psiki := \psi_i(x_i(k))$, $i=1,\mydots,n$, and $\psik = [\psi_{k,1}\; \dots \;\psi_{k,n}]^T$. Then system~\eqref{eqn:System_DT} can be rewritten (using $ \psik $ instead of $\xk$ as state variable) as:
	\begin{equation}
		\psi^{-1}(\psikp) = (I-\epsilon \Delta) \psi^{-1}(\psik) +\epsilon \pi A \psik.
		\label{eqn:System_DT_psi}
	\end{equation}	
	Let $V\!:[-1,1]^n\!\to \R_+$ be the Lyapunov function described by 
	\begin{equation*}
	V(\psik) = -\frac{1}{2} \pi \varepsilon \psik^T A \psik + \varepsilon \sum_{i} \delta_i \int_{0}^{\psiki} \psi_i^{-1}(s) \der s.
	\end{equation*}
	
	Observe that from Lemma~\ref{lemma:DT_invconvex}, $V(\psik)>0$ for all $\psik\in [-1,1]^n\setminus\{0\}$ and that $V(0)=0$. Indeed,
	\begin{align*}
		V(\psik)
		& \ge -\frac{1}{2} \pi \varepsilon \psik^T A \psik + \varepsilon \sum_{i} \frac{\delta_i}{2} \psiki^2
		\\
		& = \frac{\varepsilon}{2} \psik^T (\Delta - \pi A) \psik
		= \frac{\varepsilon}{2} \psik^T \Lpi \psik
		\ge 0.
	\end{align*}
	
	Let $\psiD = \psikp-\psik$. Computing the increment of $V$ along the trajectories gives
	\begin{align}
	V_\Delta&=V(\psikp)-V(\psik) \notag\\
	=& -\frac{\pi \varepsilon}{2} \bigl(\psikp^T A \psikp - \psik^T A \psik\bigr)
	\notag\\
	&+ \varepsilon \sum_{i} \delta_i \Bigl( \int_{0}^{\psikpi} \psi_i^{-1}(s) \der s
	 - \int_{0}^{\psiki} \psi_i^{-1}(s) \der s\Bigr)
	\notag\\
	=& -\frac{1}{2}\psiD^T (\pi \varepsilon A) \psiD - \sum_i \psiDi (\pi \varepsilon \sum_j a_{ij} \psi_{k,j})
	\notag\\& - \sum_{i} (1- \varepsilon \delta_i) \int_{\psiki}^{\psikpi}\!\!\!\!\psi_i^{-1}(s) \der s
	+ \sum_{i} \int_{\psiki}^{\psikpi}\!\!\!\!\psi_i^{-1}(s) \der s.
	\label{eqn:DT_Proof_0}
	\end{align}
	From \eqref{eqn:System_DT_psi},
	\begin{align}
	&\sum_i \psiDi (\pi \varepsilon \sum_j a_{ij} \psi_{k,j})\notag
	\\&= \sum_i \psiDi \psi_i^{-1}(\psikpi)
	- \sum_i \psiDi (1-\varepsilon \delta_i) \psi_i^{-1}(\psiki).
	\label{eqn:DT_Proof_1}
	\end{align}
	Lemma~\ref{lemma:DT_invtaylor} with $y = \psikpi$ and $y_0 = \psiki$ yields
	\begin{equation}
	\int_{\psiki}^{\psikpi} \!\!\!\psi_i^{-1}(s) \der s
	\\= \psiDi \psi_i^{-1}(\psiki) + \frac{\psiDi^2}{2}\,\der_{2}(z_i),
	\label{eqn:DT_Proof_2}
	\end{equation}
	where $\der_{2}(z_i):=\frac{1}{\psi_i'(\psi_i^{-1}(z_i))}\ge 1$ and $z_i \in [\psiki,\psikpi]$.
	While, with $y_0 = \psikpi$ and $y = \psiki$, it yields
	\begin{equation}
	\int_{\psiki}^{\psikpi} \!\!\!\psi_i^{-1}(s) \der s
	\\= \psiDi \psi_i^{-1}(\psikpi) - \frac{\psiDi^2}{2}\,\der_{2}(y_i)
	\label{eqn:DT_Proof_3}
	\end{equation}
	where $\der_{2}(y_i):=\frac{1}{\psi_i'(\psi_i^{-1}(y_i))}\ge 1$ and $y_i \in [\psiki,\psikpi]$.
	Substituting \eqref{eqn:DT_Proof_1}, \eqref{eqn:DT_Proof_2} and \eqref{eqn:DT_Proof_3} in \eqref{eqn:DT_Proof_0}, one obtains
	\begin{align*}
	&V_\Delta=
	-\frac{1}{2}\psiD^T (\pi \varepsilon A) \psiD 
	\\&\quad- \sum_i \psiDi \psi_i^{-1}(\psikpi)
		+\sum_i (1-\varepsilon \delta_i) \psiDi \psi_i^{-1}(\psiki)
	\\&\quad - \sum_{i} (1- \varepsilon \delta_i) \psiDi \psi_i^{-1}(\psiki)
		 - \sum_{i} (1- \varepsilon \delta_i) \frac{\psiDi^2}{2}\,\der_{2}(z_i)
	\\&\quad + \sum_{i} \psiDi \psi_i^{-1}(\psikpi)
	 - \sum_{i} \frac{\psiDi^2}{2}\, \der_{2}(y_i)
	 \\
	 & \le -\frac{1}{2}\psiD^T (\pi \varepsilon A) \psiD 
	 - \sum_{i} (2- \varepsilon \delta_i)\frac{\psiDi^2}{2}
	 %
	 \\
	 & = \frac{1}{2}\psiD^T \bigl (-2 I + \varepsilon \Delta - \pi \varepsilon A \bigr) \psiD
	 = \frac{1}{2}\psiD^T \bigl (-2 I + \varepsilon \Lpi \bigr) \psiD.
	\end{align*}
	The inequality holds since $1 - \varepsilon \delta_i\ge 0$, $\der_{2}(z_i)\ge 1$, and $\der_{2}(y_i)\ge 1$ for all $i$.
	Under the assumption that $-2I +\varepsilon \Lpi$ is negative definite, we obtain $V_\Delta< 0$ if $\psiD\ne 0$. Therefore the trajectories of the system~\eqref{eqn:System_DT} converge asymptotically to a fixed equilibrium point which must be the origin (see Lemma~\ref{lemma:DT_Origin}).
	\qed

\subsection{Proof of Theorem~\ref{thm:DT_summary}(ii)}
\label{proof:DT_PeriodDoubling_Pitchfork}
In this proof we follow \cite[Chapter 5]{Kuznetsov1998} and \cite{KuznetsovMeijer2005}.
The system~\eqref{eqn:System_DT} can be rewritten as
\begin{equation}
	\xkp = J_\pi \xk + F(\xk)
	\label{eqn:DT_MapProof}
\end{equation}
where $J_\pi = I -\varepsilon \Lpi$ and $F(\xk) = \pi \varepsilon A (-\xk+\psi(\xk))$.
The proof is divided into two steps. First, we assume that $\pi_1<\pid$ and prove that the system~\eqref{eqn:DT_MapProof} undergoes a pitchfork bifurcation at $\pi = \pi_1$. Then, we assume that $\pid<\pi_1$ and prove that the system~\eqref{eqn:DT_MapProof} undergoes a period-doubling bifurcation at $\pi = \pid$.

\textbf{Step 1.} Assume that $\pi_1<\pid$.
At $\pi = \pi_1$, $J_\pi$ has a simple eigenvalue at $+1$ and the corresponding eigenspace $\text{span}\{q\}$ has dimension one, where $q$ is the eigenvector of $J_{\pi_1}$ associated with $1$ (i.e., $J_{\pi_1}q = q$) normalized such that $\norm{q}_2=1$. Since $J_{\pi_1}$ is symmetric, the left eigenvector of $J_{\pi_1}$ is also $q$ (i.e., $q^T J_{\pi_1} = q^T$) .
We can decompose any vector $x\in \R^n$ as $x = u q + y$ where $u = q^T x$ and $y = x- (q^Tx) q \in (\text{span}\{q\})^\bot$.
The system of equations~\eqref{eqn:DT_MapProof} in the coordinates $(u_k,y_k)$ can be written as follows
\begin{subequations}
\begin{align}
u_{k+1} &= u_k + q^T F(u_k q+y_k)
\label{eqn:DT_MapProof_2a}
\\
y_{k+1} &=J_\pi y_k + F(u_k q+y_k) - (q^T F(u_k q+y_k)) \,q.
\label{eqn:DT_MapProof_2b}
\end{align}
\label{eqn:DT_MapProof_2}
\end{subequations}
Center manifold theory (see \cite[Chapter 5.4.2]{Kuznetsov1998}) demonstrates that the restriction of \eqref{eqn:DT_MapProof_2} to the center manifold takes the form
\begin{equation}
u_{k+1} = u_k + b \,u_k^2 + c \,u_k^3 +O(u_k^4),
\label{eqn:DT_FoldManifold}
\end{equation}
where, under the assumption~\eqref{assumption:1psiOdd}, the parameters $b$ and $c$ in \eqref{eqn:DT_FoldManifold} simplify to
\begin{equation*}
b := \frac{1}{2}q^T F_{yy}(0,\pi_1) ,\quad c :=\frac{1}{6} q^T F_{yyy}(0,\pi_1).
\end{equation*}	
Let $\ppde{\psi}{x}(x) := \diag{\ppde{\psi_1}{x_1}(x_1),\mydots,\ppde{\psi_n}{x_n}(x_n)}$ and $\ptde{\psi}{x}(x) := \diag{\ptde{\psi_1}{x_1}(x_1),\mydots,\ptde{\psi_n}{x_n}(x_n)}$. Then
\begin{equation*}
b =  \frac{1}{2}q^T F_{yy}(0,\pi_1) = \frac{\pi_1 \varepsilon}{2} q^T A \,\ppde{\psi}{x}(0) \begin{bmatrix} q_1^2 \\ \vdots \\ q_n^2\end{bmatrix}
= 0
\end{equation*}	
since $\ppde{\psi_i}{x_i}(0) =0$ for all $i$. Moreover,
\begin{align*}	
c &= \frac{1}{6} q^T F_{yyy}(0,\pi_1)
= \frac{\pi \varepsilon}{6} q^T A \ptde{\psi}{x}(0) \begin{bmatrix} q_1^3 \\ \vdots \\ q_n^3\end{bmatrix}
\notag\\
&= \frac{\varepsilon}{6} q^T \Delta \ptde{\psi}{x}(0) \begin{bmatrix} q_1^3 \\ \vdots \\ q_n^3\end{bmatrix}
= \frac{\varepsilon}{6} \sum_{i=1}^{n} \delta_i \ptde{\psi_i}{x_i}(0) q_i^4
<0,\label{eqn:DT_FoldManifold_c}
\end{align*}
since $\ptde{\psi_i}{x_i}(0)<0$ for all $i$.
Hence, \eqref{eqn:DT_FoldManifold} can be rewritten as
\begin{equation}
u_{k+1} = u_k (1-\abs{c}u_k^2) +O(u_k^4) 
\end{equation}
which means that at $\pi=\pi_1$ the system~\eqref{eqn:System_DT} undergoes a pitchfork bifurcation.

\textbf{Step 2.} Assume that $\pi_1>\pid$.
At $\pi = \pid$, $J_\pi$ has a simple eigenvalue at $-1$ and the corresponding eigenspace $\text{span}\{q\}$ has dimension one, where $q$ is the eigenvector of $J_{\pid}$ associated with $-1$ (i.e., $J_{\pid}q = q$) normalized such that $\norm{q}_2=1$.
We can decompose any vector $x\in \R^n$ as $x = u q + y$ where $u = q^T x$ and $y = x- (q^Tx) q \in (\text{span}\{q\})^\bot$.

The system of equations~\eqref{eqn:DT_MapProof} in the coordinates $(u_k,y_k)$ can be written as follows
\begin{subequations}
	\begin{align}
	u_{k+1} &= -u_k + q^T F(u_k q+y_k)
	\label{eqn:DT_MapProofPeriod_2a}
	\\
	y_{k+1} &=J_\pi y_k + F(u_k q+y_k) - (q^T F(u_k q+y_k)) \,q.
	\label{eqn:DT_MapProofPeriod_2b}
	\end{align}
	\label{eqn:DT_MapProofPeriod_2}
\end{subequations}
Center manifold theory (see \cite[Chapter 5.4.2]{Kuznetsov1998}) demonstrates that the restriction of \eqref{eqn:DT_MapProof_2} to the center manifold takes the form
\begin{equation}
u_{k+1} = -u_k + b \,u_k^2 + c \,u_k^3 +O(u_k^4),
\label{eqn:DT_PeriodManifold}
\end{equation}
where, as in \textbf{Step 1}, since each $\psi_i(\cdot)$ has odd symmetry (hence $\ppde{\psi_i}{x_i}(0) =0$ for all $i$), the parameters $b$ and $c$ in \eqref{eqn:DT_PeriodManifold} are given by
\begin{equation*}	
b =  q^T F_{yy}(0,\pid) = 0,\quad 
c = \frac{1}{6} q^T F_{yyy}(0,\pid)<0.
\end{equation*}
Hence, \eqref{eqn:DT_PeriodManifold} can be rewritten as
\begin{equation}
u_{k+1} = -u_k (1+\abs{c}u_k^2) +O(u_k^4) 
\end{equation}
which means that at $\pi=\pid$ a cycle of period 2 bifurcates from the origin.\qed

\end{document}